\documentclass[11pt]{amsart}

\usepackage{amsmath, amsthm, amssymb}
\usepackage{color}

\renewcommand{\d}{\partial}
\newcommand{\ddbar}{\sqrt{-1}\d\overline{\d}}
\newcommand{\ov}[1]{\overline{#1}}
\newcommand{\de}{\partial}
\newcommand{\db}{\overline{\partial}}
\newcommand{\Ric}{\mathrm{Ric}}
\newcommand{\mn}{\sqrt{-1}}
\newcommand{\tr}[2]{\textrm{tr}_{#1}{#2}}
\newcommand{\ti}[1]{\tilde{#1}}

\newcommand{\ve}{\varepsilon}
\numberwithin{equation}{section}

\newtheorem{thm}{Theorem}[section]
\newtheorem{prop}[thm]{Proposition}

\newtheorem{conj}[thm]{Conjecture}

\theoremstyle{definition}
\newtheorem{defn}[thm]{Definition}

\renewcommand{\leq}{\leqslant}
\renewcommand{\geq}{\geqslant}
\renewcommand{\epsilon}{\varepsilon}
\renewcommand{\le}{\leqslant}
\renewcommand{\ge}{\geqslant}
\begin{document}

\title[Gauduchon metrics with prescribed volume form]{Gauduchon metrics \\ with prescribed volume form}
\author[G. Sz\'ekelyhidi]{G\'abor Sz\'ekelyhidi}
\address{Department of Mathematics, University of Notre Dame, 255 Hurley, Notre Dame, IN 46556}
\author[V. Tosatti]{Valentino Tosatti}
\thanks{Supported in part by National Science Foundation grants DMS-1306298 and DMS-1350696 (G.Sz.), DMS-1308988 (V.T.) and DMS-1406164 (B.W.). The second-named author is supported in part by a Sloan Research Fellowship.}
\address{Department of Mathematics, Northwestern University, 2033 Sheridan Road, Evanston, IL 60208}
\author[B. Weinkove]{Ben Weinkove}
\address{Department of Mathematics, Northwestern University, 2033 Sheridan Road, Evanston, IL 60208}
\begin{abstract} We prove that on any compact complex manifold one can find Gauduchon metrics with prescribed volume form. This is equivalent to prescribing the Chern-Ricci curvature of the metrics, and thus solves a conjecture of Gauduchon from 1984.
\end{abstract}

\maketitle

\section{Introduction}

Let $M$ be a compact complex manifold of complex dimension $n$.   Suppose that $M$ admits a metric  $\alpha = \sqrt{-1} \alpha_{i\ov{j}} dz^i \wedge d\ov{z}^j>0$ which is K\"ahler (that is, $d\alpha=0$).  Yau's celebrated solution \cite{Ya} of the Calabi conjecture says that
given any smooth positive volume form $\sigma$ on $M$ with $\int_M\sigma=\int_M\alpha^n$, we can find a K\"ahler metric $\omega$ with this prescribed volume form
\begin{equation}
\label{yau}\omega^n=\sigma.
\end{equation}
Moreover, there exists such a metric so that $[\omega]=[\alpha]$ in $H^2(M,\mathbb{R})$, and with this cohomological constraint the metric $\omega$ is unique.

Furthermore, Yau's Theorem is equivalent to a statement about the first Chern class $c_1(M)$.  Namely, given any smooth representative $\Psi$ of $c_1(M)$, there exists a unique K\"ahler metric $\omega$ cohomologous to $\alpha$ such that
\begin{equation}\label{yaur}\Ric(\omega)=\Psi,
\end{equation}
where $\Ric(\omega)$ is the Ricci form of the K\"ahler metric $\omega$.
Indeed, this follows immediately from the definition of $c_1(M)$ and by applying the operator $-\ddbar \log$ to (\ref{yau}).

It is natural to investigate whether similar results hold when $M$ does not admit a K\"ahler metric, but only a Hermitian metric $\alpha$.   If we do not impose any constraint on the class of Hermitian metrics that we consider, then \eqref{yau} can be trivially solved by a conformal change of metric. However, there is a natural class of Hermitian metrics which exist on all compact complex manifolds, namely {\em Gauduchon metrics}. A Hermitian metric $\alpha$ is called Gauduchon if \[ \de\db(\alpha^{n-1})=0, \] and a classical result of Gauduchon \cite{Ga0} says that every Hermitian metric is conformal to a Gauduchon metric (uniquely up to scaling, when $n\geq 2$). In particular, if we restrict our attention to Gauduchon metrics then we cannot use nontrivial conformal changes.

Motivated by Yau's Theorem, in 1984 Gauduchon \cite[IV.5]{Ga} posed the following conjecture:

\begin{conj}\label{conj1}
Let $M$ be a compact complex manifold and $\Psi$ a closed real $(1,1)$ form on $M$ with $[\Psi]=c_1^{\rm{BC}}(M)\in H^{1,1}_{\mathrm{BC}}(M,\mathbb{R})$. Then there is a Gauduchon metric
$\omega$ on $M$ with
\begin{equation}\label{todo}\Ric(\omega)=\Psi.\end{equation}
\end{conj}

To explain our notation here,
$$H^{1,1}_{\mathrm{BC}}(M,\mathbb{R})=\frac{\{d{\textrm-}{\rm closed\ real\ }(1,1){\rm \ forms}\}}{\{\ddbar\psi\ |\ \psi\in C^\infty(M,\mathbb{R})\}},$$
denotes the (finite dimensional) {\em Bott-Chern cohomology group}, and $\Ric(\omega)$ is the Chern-Ricci form of $\omega$,  which is locally given by
$$\textrm{Ric}(\omega)=-\ddbar\log \det g,$$ where we write $\omega = \sqrt{-1} g_{i\ov{j}}dz^i \wedge d\ov{z}^j$.
 It is a closed real $(1,1)$ form and its first Bott-Chern cohomology class
$c_1^{\rm{BC}}(M)=[\Ric(\omega)]\in H^{1,1}_{\mathrm{BC}}(M,\mathbb{R})$ is immediately seen to be independent of the choice of $\omega$.

In the spirit of Yau's Theorem, we restate Conjecture \ref{conj1} as an equivalent statement about the existence of Gauduchon metrics with prescribed volume form:

\begin{conj}\label{conj2}
Let $M$ be a compact complex manifold and $\sigma$ a smooth positive volume form. Then there is a Gauduchon metric $\omega$ on $M$ with
\begin{equation}
\label{todo2}\omega^n=\sigma.\end{equation}
\end{conj}

The equivalence with Conjecture \ref{conj1} follows  by applying the operator $-\ddbar\log$ to (\ref{todo2}).

Our result, Theorem \ref{conj3} below, gives a proof of Conjecture \ref{conj1} (and hence also Conjecture \ref{conj2}).  Moreover, our result strengthens the conjecture by imposing a cohomological constraint on the solution $\omega$.  Before we state our results, we make some remarks about Conjecture \ref{conj1}:\\

\noindent (1) When $M$ is K\"ahler this conjecture follows from Yau's Theorem.\\

\noindent (2) When $n=2$ the conjecture was proved by Cherrier \cite{Ch} in 1987 by solving a complex Monge-Amp\`ere equation (see also \cite{TW1, GL} for different proofs).\\

\noindent (3) More recently the second and third-named authors \cite{TW3} proved Conjecture \ref{conj1} when $M$ admits an {\em astheno-K\"ahler} metric, i.e. a Hermitian metric $\alpha$ with $\de\db(\alpha^{n-2})=0$ (a condition introduced in \cite{JY}).\\

\noindent (4) Clearly there can be no uniqueness in Conjecture \ref{conj1} as stated.\\

\noindent (5) In \cite{TW4} the second and third-named authors proved that given a Hermitian metric $\alpha$ one can always find another Hermitian metric $\omega$ of the form $\omega=\alpha+\ddbar u$ for $u \in C^\infty(M,\mathbb{R})$, solving \eqref{todo}. If $n=2$ then $\alpha$ Gauduchon implies that $\omega$ is also Gauduchon (and this equation was solved in \cite{Ch}), but this is no longer the case when $n\geq 3$.  Hence the result of \cite{TW4} does not help to solve Conjecture \ref{conj1} in dimension $3$ or higher. \\

\noindent (6) A consequence of Conjecture \ref{conj1} is that $c_1^{\rm{BC}}(M)=0$ holds if and only if there exist Chern-Ricci-flat Gauduchon metrics on $M$. More information about these ``non-K\"ahler Calabi-Yau'' manifolds
can be found in \cite{To}.\\

We now state our main results.  We first introduce some terminology concerning cohomology classes of $(n-1,n-1)$ forms.  Define the \emph{Aeppli cohomology group}
$$H^{n-1,n-1}_{\mathrm{A}}(M,\mathbb{R})=\frac{\{\de\db{\textrm-}{\rm closed\ real\ }(n-1,n-1){\rm \ forms}\}}{\{\de\gamma+\ov{\de\gamma}\ |\ \gamma\in \Lambda^{n-2,n-1}(M)\}}.$$
This space is naturally in duality with the Bott-Chern cohomology group we considered earlier, with the nondegenerate pairing $H^{n-1,n-1}_{\mathrm{A}}(M,\mathbb{R})\otimes H^{1,1}_{\mathrm{BC}}(M,\mathbb{R})\to\mathbb{R}$ given by wedge product and integration over $M$ (see e.g. \cite{An}). If $\alpha_0$ is a Gauduchon metric then $\alpha_0^{n-1}$ defines a class $[\alpha_0^{n-1}]\in H^{n-1,n-1}_{\mathrm{A}}(M,\mathbb{R})$.

We prove:
\begin{thm}\label{conj3}
Let $M$ be a compact complex manifold with a Gauduchon metric $\alpha_0$, and $\Psi$ a closed real $(1,1)$ form on $M$ with $[\Psi]=c_1^{\rm{BC}}(M)\in H^{1,1}_{\mathrm{BC}}(M,\mathbb{R})$.
Then there exists a Gauduchon metric $\omega$ satisfying $[\omega^{n-1}]=[\alpha_0^{n-1}]$ in $H^{n-1,n-1}_{\mathrm{A}}(M,\mathbb{R})$ and
\begin{equation}\label{ricc}
\Ric(\omega)=\Psi.
\end{equation}
\end{thm}

This result immediately implies Conjectures \ref{conj1} and \ref{conj2}.

 In \cite{TW3}, the second and third-named authors observed that to solve Theorem \ref{conj3}  it is enough to solve a certain partial differential equation, which was also independently introduced by Popovici \cite{Po}.  This equation is a variant of one introduced by Fu-Wang-Wu \cite{FWW1} and related to Harvey-Lawson's notion of $(n-1)$-plurisubharmonic functions \cite{HL1, HL2}.

Namely, we seek a Hermitian metric $\omega$ on $M$ with the property that
$$\omega^{n-1}=\alpha_0^{n-1}+\de\gamma+\ov{\de\gamma},$$
where
$$\gamma=\frac{\mn}{2}\db u \wedge\alpha^{n-2},$$
$u \in C^\infty(M,\mathbb{R})$ and $\alpha$ is a background Gauduchon metric. Clearly, by construction, the metric $\omega$ is Gauduchon assuming $\alpha_0$ is Gauduchon. Substituting, we see that
\begin{equation}\label{eqa}
\omega^{n-1}=\alpha_0^{n-1}+\ddbar u \wedge\alpha^{n-2}+\mathrm{Re}\left(\mn \de u \wedge\db(\alpha^{n-2})\right),
\end{equation}
while if we write
$$\Ric(\omega)=\Psi+\ddbar F,$$
then \eqref{ricc} is equivalent to
$$\omega^n=e^{F+b}\alpha^n,$$
for some constant $b\in\mathbb{R}$. This is exactly the equation that we solve, thus resolving \cite[Conjecture 1.5]{TW3} and \cite[Question 1.2]{Po}.

\begin{thm}\label{main2}
Let $M$ be a compact complex manifold with $\dim_{\mathbb{C}}M=n\geq 2$, equipped with a Hermitian metric $\alpha_0$ and a Gauduchon metric $\alpha$. Given a smooth function $F$ on $M$ we can find a unique $u \in C^\infty(M,\mathbb{R})$ with $\sup_M u =0$, and a unique $b\in\mathbb{R}$ such that the Hermitian metric $\omega$ defined by
$$\omega^{n-1}:=\alpha_0^{n-1}+\ddbar u \wedge\alpha^{n-2}+\mathrm{Re}\left(\mn \de u \wedge\db(\alpha^{n-2})\right)>0,$$
satisfies
\begin{equation}\label{cma}
\omega^n=e^{F+b}\alpha^n.
\end{equation}
\end{thm}

Clearly, as we just described, Theorem \ref{conj3} follows from this result if we take $\alpha_0=\alpha$ Gauduchon.  We make some remarks about Theorem \ref{main2}.\\

\noindent
(1) In the case when $\alpha$ is K\"ahler, or more generally if the linear  term involving $\partial u$ is removed, the equation (\ref{cma}) reduces to the Monge-Amp\`ere equation for $(n-1)$-plurisubharmonic functions, solved by the second and third-named authors \cite{TW2, TW3} (see also \cite{FWW2} for earlier partial results).  \\

\noindent
(2) In the case when $n=2$ this equation reduces to the complex Monge-Amp\`ere equation, solved in \cite{Ch} (see also \cite{TW1}).\\

\noindent
(3) It was shown in \cite{TW3} that Theorem \ref{main2} can be reduced to a second order a priori estimate of the form (cf. \cite{HMW})
$$\sup_M |\sqrt{-1} \partial \ov{\partial} u|_{\alpha} \le C(1+ \sup_M| \nabla u|^2_{\alpha}),$$
for solutions $u$ of (\ref{cma}).  This is precisely the estimate we prove in this paper.\\

\noindent
(4) If in Theorem \ref{main2} we assume that $\alpha_0$ is \emph{strongly Gauduchon} in the sense of Popovici \cite{P1}, namely that $\ov{\partial} (\alpha_0^{n-1})$ is $\partial$-exact, then by construction so is the solution $\omega$.  Thus we also get a Calabi-Yau-type theorem for strongly Gauduchon metrics. More applications of this theorem can be found in \cite{Po}.\\

\noindent
(5) Our method of proof of Theorem \ref{main2} can also be used to solve an equation introduced by Fu-Wang-Wu \cite{FWW1} in certain cases. Suppose we have a compact Hermitian manifold $(M,\alpha_0)$ and we seek a Hermitian metric $\omega$ solving \eqref{cma} with the property that
$$\omega^{n-1}=\alpha_0^{n-1}+\ddbar(u\, \alpha^{n-2}),$$
for some $u \in C^\infty(M,\mathbb{R})$ and some  Hermitian metric $\alpha$. This setup is particularly interesting because if $\alpha_0$ is balanced (i.e. $d(\alpha_0^{n-1})=0$, see \cite{Mi}), then so is $\omega$, and one obtains a Calabi-Yau theorem for balanced metrics (see also \cite[Section 4]{To}). When $\alpha$ is K\"ahler this setup reduces to the setting of item (1).  If we instead assume that $\alpha$ is astheno-K\"ahler, then we see that
\begin{equation}\label{eqb}
\omega^{n-1}=\alpha_0^{n-1}+\ddbar u \wedge\alpha ^{n-2}+2\mathrm{Re}\left(\mn \de u \wedge\db(\alpha^{n-2})\right),
\end{equation}
which differs from \eqref{eqa} just for a factor of $2$.
Therefore this problem falls into our general framework (see Theorem \ref{main} below), and we conclude that we have uniform {\em a priori} estimates for solutions of this equation.
The exact same argument as in \cite[Theorem 1.7]{TW3} using the continuity method then shows that the equation is indeed solvable. In the case when we choose $\alpha_0$  to be balanced, this gives a proof of \cite[Conjectures 4.1 and 4.2]{To} assuming that $M$ admits astheno-K\"ahler metrics.
However, we should remark that we are not aware of any example of a non-K\"ahler compact complex manifold which admits both balanced and astheno-K\"ahler metrics.\footnote{After this paper was posted, and prompted by our remark, explicit examples were constructed in \cite{FGV,LU} in all complex dimensions $\geq 4$.} \\

\noindent
(6) The same argument as the proof of Theorem~\ref{main2}
 also allows us to find a Gauduchon metric $\omega$ solving the ``complex-Hessian'' equation
$$\omega^k\wedge\alpha^{n-k}=e^{F+b}\alpha^n,$$
for any $1\leq k\leq n$, see also \cite[Proposition 24]{Sz} for the case of $(n-1)$-plurisubharmonic functions, and \cite{DK, HMW} for the standard K\"ahler case where $\omega=\alpha+\ddbar u$.\\

\noindent
(7) The complex setting is very different from the real analogue of (\ref{cma}), treated for example in more generality in \cite{GuJ}.  The underlying reason is the two different types of complex derivatives.  In our case the special structure of the gradient term in (\ref{cma}) plays a key role.\\

In fact, Theorem \ref{main2} follows from a much more general result
where we consider a large class of fully non-linear second-order elliptic equations on Hermitian
manifolds. This result is analogous to the main result in \cite{Sz},
giving a priori estimates in the presence of a suitable
subsolution. We will state this as Theorem~\ref{main} in
Section~\ref{sec:bg}. This result fits into a large body of work on
fully non-linear second order elliptic equations, going back to the
work of Caffarelli-Nirenberg-Spruck~\cite{CNS} on the Dirichlet
problem on domains in $\mathbb{R}^n$. Some other works on this
topic include \cite{Chu, CSz, Gu, Gu2, GuJ, GL, Li, Li2, Li3, Ni, PSS, Su1, Su2, Tru, Urb, ZhD, Zh, ZZ}.

In our proof of Theorem \ref{main} we  use  some of the language and approaches of the recent paper of the first-named author \cite{Sz}.  However, if one is only interested in a direct proof of Theorem \ref{main2}, one can equally well use the language of \cite{TW3}.  In any case, the key new ingredient is an understanding of the structure of the term  $\mathrm{Re}\left(\mn \de u \wedge\db(\alpha^{n-2})\right)$.\\

The paper is organized as follows. In section \ref{sec:bg} we will introduce some notation and state our main technical theorem \ref{main}. The proof of this theorem will be given in section \ref{proof}, and in section \ref{sec:last} we show how this implies Theorem  \ref{main2}.\\

 As the present work neared completion, we were informed that Bo Guan and Xiaolan Nie have a work in progress on related results.\\

\noindent
{\bf Acknowledgments. }We thank the referee for useful remarks.

\section{Background and the General Setting}\label{sec:bg}

Let  $(M,\alpha)$ be a compact Hermitian manifold of complex dimension $n$ and write
$$\alpha = \sqrt{-1} \alpha_{i\ov{j}} dz^i \wedge d\ov{z}^{j}>0.$$
Fix a background $(1,1)$ form $\chi = \sqrt{-1} \chi_{i\ov{j}}dz^i \wedge d\ov{z}^j$ which is not necessarily positive definite.  Let $W_{i\bar j}(\nabla u)$ be a Hermitian tensor which depends linearly on $\nabla u$. For $u: M\to\mathbb{R}$
define a new tensor $g_{i\ov{j}}$ by
\begin{equation}
 \label{eq:u4} g_{i\bar j} := \chi_{i\bar j} + u_{i\bar j} + W_{i\bar j}. \end{equation}
Note that we do not assume that $(g_{i\ov{j}})$ is positive definite.
We will study equations for $g$, where $W$ has a special structure
related to the equation \eqref{cma}. To define this let us write
\begin{equation}
 \label{eq:u5} \tilde{g}_{i\bar j} := P_\alpha(g_{i\bar j}) = \frac{1}{n-1}\left( (\mathrm{tr}_\alpha g)
  \alpha_{i\bar j} - g_{i\bar j}\right), \end{equation}
where $P_\alpha$ is an operator on tensors, depending on the fixed metric $\alpha$, defined by the second equality in (\ref{eq:u5}).  As an aside, if $\alpha$ is the Euclidean metric on $\mathbb{C}^n$ then the condition $P_{\alpha}(u_{i\ov{j}}) \ge 0$ is equivalent to saying that $u$ is $(n-1)$-plurisubharmonic, in the sense of Harvey-Lawson \cite{HL1}.

Observe that, writing $\Delta = \alpha^{k\ov{\ell}} \partial_k \partial_{\ov{\ell}}$,
\begin{equation}
 \label{eq:u6}\tilde{g}_{i\bar j} = \tilde{\chi}_{i\bar j} + \frac{1}{n-1}\left(
  (\Delta u)\alpha_{i\bar j} - u_{i\bar j}\right) + Z_{i\bar j} \end{equation}
for $Z$ given by
\begin{equation}
 \label{eq:c1} Z_{i\bar j} := P_\alpha(W_{i\ov{j}})= \frac{1}{n-1}\left( (\mathrm{tr}_\alpha
  W)\alpha_{i\bar j} - W_{i\bar j} \right), \end{equation}
  and similarly, $\tilde{\chi}_{i\ov{j}} = P_\alpha(\chi_{i\ov{j}})$.
  Note that we can also write $W$ explicitly in terms of $Z$
 \begin{equation}
 \label{eq:cz} W_{i\ov{j}} = (\mathrm{tr}_{\alpha} Z) \alpha_{i\ov{j}} - (n-1) Z_{i\ov{j}}.
 \end{equation}
A crucial assumption we make is that $W$ depends on $\nabla
u$ in the following way:
we assume that the tensor $Z$ has the form
\begin{equation}\label{ZZZ}
Z_{i\ov{j}}=Z_{i\ov{j}}^p u_p +\ov{Z_{j\ov{i}}^p u_p},
\end{equation}
for some tensor $Z_{i\ov{j}}^p$, independent of $u$. In addition we have:

\bigskip
\noindent {\bf Assumption for $W$:} In orthonormal coordinates for $\alpha$ at any given point, the component
$Z_{i\bar j}$ is independent of $u_{\bar i}$ and $u_j$ (in other words $Z_{i\ov{j}}^j=0$ for all $i,j$), and
$\nabla_{i}Z_{i\ov{i}}$ is independent of $u_{\ov{i}}$ (in other words $\nabla_{\ov{i}} Z_{i\ov{i}}^i=0$ for all $i$).
\bigskip

Here $\nabla$ is the Chern connection of $\alpha$. This assumption expresses a certain skew-symmetry requirement for the tensor $W$.
This assumption is satisfied for  the $(n-1)$-plurisubharmonic Monge-Amp\`ere equation, the case of most interest to us, see \eqref{eq:Zdefn2} below,
the key reason being that the torsion tensor is skew-symmetric.

Let us record here a few consequences of this assumption, which will be used later. Taking $\nabla_p$ of \eqref{ZZZ}, setting $i=j$, evaluating
at that point and using that $Z_{i\ov{i}}^i=0$ we see that $\nabla_p Z_{i\ov{i}}$ is independent of $u_{ii}$ and $u_{i\ov{i}}$ (at that point, in orthonormal coordinates for $\alpha$).  Here the subscripts of $u$ denote ordinary partial derivatives.
Similarly, $\nabla_i Z_{p\ov{i}}$ is independent of $u_{ii}$. Taking two covariant derivatives we have
\[\begin{split}\nabla_{\ov{i}} \nabla_i Z_{i\ov{i}}=&\nabla_{\ov{i}} \nabla_i Z_{i\ov{i}}^p u_p +Z_{i\ov{i}}^p \nabla_{\ov{i}} \nabla_i u_p+\nabla_i Z_{i\ov{i}}^p \nabla_{\ov{i}}u_p+
\nabla_{\ov{i}}Z_{i\ov{i}}^p \nabla_i u_p\\
&+\nabla_{\ov{i}} \nabla_i\ov{Z_{i\ov{i}}^p} u_{\ov{p}}+\ov{Z_{i\ov{i}}^p} \nabla_{\ov{i}} \nabla_i u_{\ov{p}}+\nabla_i \ov{Z_{i\ov{i}}^p}\nabla_{\ov{i}} u_{\ov{p}}+
\nabla_{\ov{i}}\ov{Z_{i\ov{i}}^p}\nabla_i u_{\ov{p}},
\end{split}\]
and evaluating at our point and using the assumptions $\nabla_{\ov{i}} Z_{i\ov{i}}^i=0$ and $Z_{i\ov{i}}^i=0$, we see that $\nabla_{\ov{i}} \nabla_i Z_{i\ov{i}}$ is independent of
$u_{ii}, u_{\ov{i}\ov{i}}, u_{ii\ov{i}}$ and $u_{i\ov{i}\ov{i}}$.

Given a smooth function $h$, we study equations of the form
\[ F(A)=h,\]
where $A$ is the endomorphism
$A^i_j = \alpha^{i\bar p}g_{j\bar p}$ of the holomorphic tangent
bundle, which is Hermitian with respect to the inner product defined by $\alpha$, and $F(A)$ is a symmetric function of the eigenvalues
$\lambda_1,\ldots,\lambda_n$ of $A$:
\begin{equation}
\label{use} F(A) = f(\lambda_1,\ldots,\lambda_n). \end{equation}
We assume that our operator $F$ has the special form $F(M) =
\widetilde{F}(P(M))$, where $P(M) = (n-1)^{-1}\big[ \mathrm{Tr}(M)I -
M\big]$, analogous to $P_{\alpha}$ above, and
\[\widetilde{F}(B)=\ti{f}(\mu_1,\dots,\mu_n),\]
where $\mu_1,\dots,\mu_n$ are the eigenvalues of $B$, and $\ti{f}$ is another symmetric function.
In terms of eigenvalues, this means that
\begin{equation} \label{eq:u2} f(\lambda_1,\ldots, \lambda_n) = (\tilde{f} \circ P) (\lambda_1, \ldots, \lambda_n),
\end{equation}
where we are writing $P$ for the map $\mathbb{R}^n \rightarrow \mathbb{R}^n$ induced on diagonal  matrices by the matrix map $P$ above.
Explicitly, writing $\mu = P(\lambda)$ for $\lambda,\mu\in\mathbb{R}^n$ the corresponding $n$-tuples, we have
\begin{equation}
 \label{eq:u1} f(\lambda_1, \ldots, \lambda_n) = \tilde{f}(\mu_1, \ldots, \mu_n), \quad \textrm{for }\mu_k = \frac{1}{n-1}\sum_{i\ne k} \lambda_i. \end{equation}

\bigskip
\noindent{\bf Assumptions for $\widetilde{f}$ and $h$:} we
make the following assumptions on $\widetilde{f}$, and the function
$h$ in our equation:
\begin{itemize}
\item[(i)] $\widetilde{f}$ is defined on an open symmetric convex cone
  $\widetilde{\Gamma}\subsetneq \mathbb{R}^n$, containing the positive orthant
  $\Gamma_n = \{ (x_1, \ldots, x_n) \in \mathbb{R}^n \ | \ x_i>0, \ i=1, \ldots, n \}.$
\item[(ii)] $\widetilde{f}$ is symmetric, smooth, concave, and increasing,
  i.e. its partials satisfy $\widetilde{f}_i > 0$ for all $i$.
\item[(iii)] $\sup_{\d\widetilde{\Gamma}} \widetilde{f} < \inf_M h$.
\item[(iv)] For all $\mu\in \widetilde{\Gamma}$ we have
  $\lim_{t\to\infty} \widetilde{f}(t\mu) = \sup_{\widetilde{\Gamma}}
  \widetilde{f}$, where both sides are allowed to be $\infty$.
  \item[(v)] $h$ is a smooth function on $M$.
\end{itemize}

Define the cone $\Gamma\subset \mathbb{R}^n$ by $\Gamma =
P^{-1}(\widetilde{\Gamma})$.  Observe that $P$ maps $\Gamma_n$ into $\Gamma_n$.
It is then easy to see that the function $f = \tilde{f} \circ P:\Gamma\to\mathbb{R}$
 satisfies exactly the same conditions as
$\widetilde{f}$. In particular some of the results of \cite{Sz} can be applied
to the equation $F(A) = h$. We need the following definition (see
Remark 8 in \cite{Sz} to see the equivalence with the definition
there), which is a modification of a notion introduced by Guan \cite{Gu2}.

\begin{defn} We say that $u$ is a $\mathcal{C}$-subsolution for the
  equation $F(A) = h$ if the following holds. Let $g_{i\bar j}$ be
  defined as in \eqref{eq:u4}. We require that for every point $x\in
  M$, if $\lambda = (\lambda_1,\ldots, \lambda_n)$ denote the eigenvalues of the
  endomorphism $\alpha^{i\bar p}g_{j\bar p}$ at $x$, then for all $i=1,\dots,n$ we have
  \[ \lim_{t\to\infty} f(\lambda + t\mathbf{e}_i) > h(x). \]
  Here $\mathbf{e}_i$ denotes the $i^\mathrm{th}$ standard basis
  vector. Note that part of the requirement is that $\lambda +
  t\mathbf{e}_i \in \Gamma$ for sufficiently large $t$, for the limit
  to be defined.
\end{defn}

With this background, our main estimate is the following, analogous to
the main result in
\cite{Sz}. We will give the proof in section~\ref{proof}.

\pagebreak[3]
\begin{thm}\label{main}
   Suppose that $\underline{u}$ is a $\mathcal{C}$-subsolution for the
   equation $F(A)=h$, and $u$ is a smooth solution, normalized by
   $\sup_Mu=0$. Suppose that $F$ and $h$ satisfy the assumptions
   above, including the assumption for the gradient term $W$.
   Then for each $k=0,1,2,\ldots$, we have an
   estimate $\Vert u \Vert_{C^{k}(M,\alpha)} \leq C_k$, with
   constant $C_k$ depending on $k$, on the background data $M,\alpha,
   \chi, F, h$,  the
   coefficients of $W$ and the subsolution $\underline{u}$.
\end{thm}

The case of primary interest for us is   equation \eqref{cma}, which corresponds to the symmetric function
\begin{equation}
\label{eq:lm}
\widetilde{f}(\mu_1,\ldots,\mu_n) = \log (\mu_1\cdot\ldots \cdot
\mu_n),\end{equation}
on the positive orthant $\widetilde{\Gamma}=\Gamma_n$.
 It is straightforward to check that $\widetilde{f}$ satisfies
the conditions above. Indeed, $\widetilde{f}$ converges to
$-\infty$ on the boundary $\d\Gamma_n$, so (iii) is satisfied, and for
(iv) it is enough to note that $\widetilde{f}(t\mu) =
\widetilde{f}(\mu) + n\log t$, which converges to $\infty$ as
$t\to\infty$.

In addition, if $\mu\in \Gamma_n$, then we also have
\[ \lim_{t\to\infty} \widetilde{f}(\mu + t\mathbf{e}_i) = \infty, \]
for all $i$. This means that for a function $u$ to be a
$\mathcal{C}$-subsolution for this equation, the only requirement is
that at each point the eigenvalues $\lambda$ of $\alpha^{i\bar
  p}g_{j\bar p}$ satisfy $P(\lambda)\in \Gamma_n$. In other words, the
requirement is that $\widetilde{g}_{i\bar j}$, defined in
\eqref{eq:u5}, is positive definite.

Note that if $\underline{u}$ is a $\mathcal{C}$-subsolution, then
replacing $\chi$ by
\[ \chi'_{i\bar j} = \chi_{i\bar j} + \underline{u}_{i\bar j} +
W_{i\bar j}(\nabla \underline{u}), \]
we can assume that $\underline{u}=0$. The important consequence of $0$
being a $\mathcal{C}$-subsolution is the following, which follows from
 Proposition 6 and Lemma 9 in \cite{Sz}.

\begin{prop}\label{prop:subsol}
Suppose that $0$ is a $\mathcal{C}$-subsolution for the equation
$F(A)=h$, and $u$ is a solution. Define $g_{i\bar j}$ as in
\eqref{eq:u4}.
There are constants $R, \kappa > 0$, independent of $u$, with the following
property. Let $x\in M$, and choose orthonormal coordinates for
$\alpha$ at $x$, such that $g$ is diagonal, with eigenvalues
$\lambda = (\lambda_1, \ldots,  \lambda_n)$. If $|\lambda| > R$, then there
are two possibilities:
\begin{enumerate}
\item[(a)] We have
\[ \sum_k f_k(\lambda)\, [\chi_{k\bar k} - \lambda_k] > \kappa\sum_k
f_k(\lambda). \]
\item[(b)] Or, $f_k(\lambda) > \kappa \sum_i f_i(\lambda)$ for all $k$.
\end{enumerate}
In addition $\sum_k f_k(\lambda) > \kappa$.
\end{prop}

We collect some other basic properties of the functions $f$ and
$\widetilde{f}$. Suppose that $\lambda \in \Gamma$ with
 $\lambda_1 \geq\ldots \geq \lambda_n$. Then
$\mu_1\leq\ldots\leq \mu_n$, and so by property (ii), $\widetilde{f}_1 \geq \ldots \geq
\widetilde{f}_n>0$ (see e.g. \cite[p.12]{Sz}). We have
\begin{equation} \label{fk}
f_k = \frac{1}{n-1}\sum_{i\ne k} \widetilde{f}_i,
\end{equation}
which implies that $0<f_1  \le \cdots \le f_n$.  Also, for $k > 1$,
\begin{equation}
 \label{eq:fk1}
0<\frac{\widetilde{f}_1}{n-1} \leq f_k \leq \widetilde{f}_1, \end{equation}
i.e. the $f_k$ for $k > 1$ are all comparable, while $f_1$ may be
relatively small. In addition, from (\ref{fk}) with $k=1$, we obtain
 \begin{equation}
  \label{eq:fk2}
 \widetilde{f}_i \leq (n-1)f_1, \textrm{ for } i
> 1.
\end{equation}

Proposition~\ref{prop:subsol} is easy to verify directly in the case of  equation \eqref{cma}, where
\[ f(\lambda) = \log (\mu_1\cdot\ldots\cdot \mu_n), \]
with $\mu_k$ defined as in \eqref{eq:u1}. Indeed, in this case
\[ \widetilde{f}_i(\mu) = \frac{1}{\mu_i}, \]
and
\[ f_k(\lambda) = \frac{1}{n-1}\sum_{i\ne k} \frac{1}{\mu_i}. \]

The function $0$ being a $\mathcal{C}$-subsolution means that
$\widetilde{\chi}$ in \eqref{eq:u6} is positive definite. We have
\[ \sum_k f_k(\lambda) \chi_{k\bar k} = \sum_i \frac{1}{\mu_i}
\widetilde{\chi}_{i\bar i} > \tau\sum_i \frac{1}{\mu_i} =
\tau\sum_k f_k(\lambda), \]
for some $\tau > 0$ depending on a lower bound for
$\widetilde{\chi}$. We also have
\[ \sum_k \lambda_k f_k(\lambda) = n. \]
It follows that we have the alternative (a) in
Proposition~\ref{prop:subsol} whenever there is one sufficiently small
$\mu_i$, which by the equation $f(\lambda) = h$ is equivalent to
having at least one large $\mu_i$, i.e. at least one large
$\lambda_i$. In addition
\[ \begin{aligned}\sum_k f_k(\lambda) &= \sum_{i=1}^n \frac{1}{\mu_i} \\
&\geqslant n (\mu_1\cdot\ldots\cdot \mu_n)^{-1/n} \\
&= ne^{-h / n}, \end{aligned}\]
so that the final claim in Proposition~\ref{prop:subsol} also holds.

\section{Proof of the main estimate}\label{proof}

In this section we give the proof of Theorem \ref{main}.

First of all, a uniform bound $\|u\|_{L^\infty(M)}\leq C$ can be
obtained by a simple modification of the argument
in \cite[Proposition 10, Remark 12]{Sz}, which is itself inspired by B\l ocki's proof of the $L^\infty$ estimate in Yau's Theorem \cite{Bl}.
In the setting of equation (\ref{cma}), the $L^{\infty}$ estimate of $u$ was first proved in \cite{TW3}, using a different method more analogous to the arguments in \cite{Ya, Ch, TW4, TW2}.

Our main goal is the following estimate:
\begin{equation}
\label{mainbd}\sup_M |\ddbar u|_{\alpha}\leq C(\sup_M|\nabla u|^2_\alpha+1),\end{equation}
for a constant $C$ depending only on the fixed data of Theorem \ref{main}. We remark that an estimate of this form was proved in the context of the complex Hessian equations by Hou-Ma-Wu \cite{HMW}, making use of ideas of Chou-Wang \cite{CW}.  For the $(n-1)$-plurisubharmonic equation (namely, equation \eqref{cma} without the linear term in $\partial u$),
an estimate of this type was proved in \cite{TW2, TW3}.  This was then generalized much further in
\cite{Sz}, where the estimate was shown to hold for a large class of equations.  Our proof begins along similar lines to  these papers.  The new difficulty comes from the linear term in $\partial u$, which, fortunately, has a special structure that we can exploit.

 In fact, the estimate (\ref{mainbd}) is equivalent to the bound
\[ \lambda_1 \le C K, \]
where  $K = 1 + \sup_M
|\nabla u|_{\alpha}^2$  and $\lambda_1$ is the largest eigenvalue of $A= (A^i_j)=(\alpha^{i\bar p}g_{j\bar p})$.  Indeed, our assumption on the cone $\Gamma$ implies that
$\sum_i \lambda_i>0$ (see Caffarelli-Nirenberg-Spruck~\cite{CNS}).  Then if $\lambda_1$ is bounded from above by $CK$ then so is $|\lambda_i|$ for all $i$, giving the same bound for $\sup_M |\ddbar u|_{\alpha}$.

We consider the function
\[ H = \log \lambda_1 + \phi(|\nabla u|^2_\alpha) +
\psi(u), \]
where $\phi$ is defined by
\[ \phi(t) = -\frac{1}{2}\log\left( 1 - \frac{t}{2K}\right) \]
so that $\phi(| \nabla u|^2_{\alpha}) \in [0, (\log 2)/2] $ satisfies
\[ (4K)^{-1} < \phi' < (2K)^{-1}, \quad \phi'' = 2\phi'^{2} > 0,   \]
and $\psi$ is defined by
\begin{equation}
 \label{eq:psid}
\psi(t) = D_1e^{-D_2t}, \end{equation}
for sufficiently large uniform constants $D_1,D_2 > 0$ to be chosen later.  By the $L^{\infty}$ bound on $u$, the quantity $\psi(u)$ is uniformly bounded.

We remark that we follow \cite{Sz} by computing with the largest eigenvalue $\lambda_1$ instead of the analogous quantity in \cite{TW3}, but in fact either quantity works, at least in the case of equation \eqref{cma}. Also, note that while the function $\phi$ here coincides with that in \cite{HMW} (and also in \cite{TW2,TW3}), our choice of $\psi$ is crucially different.

We work at a point  where $H$ achieves its maximum, in orthonormal complex
coordinates for $\alpha$ centered at this point, such that $g$ is diagonal and $\lambda_1 =
g_{1\bar 1}$.
The quantity $H$ need not be smooth at this maximum point because the largest eigenvalue of $A$
may have eigenspace of dimension larger than $1$. To take care of this, we carry out a perturbation argument as in \cite{Sz},  choosing local
coordinates such that $H$ achieves its maximum at the origin, where $A$ is diagonal with eigenvalues $\lambda_1 \ge \cdots \ge \lambda_n$, as before.
We fix a diagonal matrix $B$ with $B^1_1=0$ and $0<B^2_2<\dots<B^n_n$,
and we define $\ti{A}=A-B$, and denote its eigenvalues by
$\tilde{\lambda}_1, \ldots, \tilde{\lambda}_n$.

At the origin we have
$$\ti{\lambda}_1=\lambda_1,\quad \ti{\lambda}_i=\lambda_i-B^i_i, \quad i>1,$$
and $\tilde{\lambda}_1 > \tilde{\lambda_2} > \cdots > \tilde{\lambda}_n$. As discussed above, our assumption on the cone $\Gamma$ implies that
$\sum_i \lambda_i>0$
and we fix the matrix $B$ small enough so that
\[ \sum_i \ti{\lambda}_i>-1. \]
We can choose such $B$ such that, in addition,
\begin{equation}\label{bdd}
\sum_{p>1}\frac{1}{\lambda_1-\ti{\lambda}_p}\leq C,
\end{equation}
for some fixed constant $C$ depending on the dimension $n$.
Now, after possibly shrinking the chart, the quantity
\[ \ti{H} = \log \ti{\lambda}_1 + \phi(|\nabla u|^2_\alpha) +
\psi(u), \]
is smooth on the chart and achieves its maximum at the origin.  We will apply the maximum principle to $\tilde{H}$.  Our goal is to obtain the bound $\tilde{\lambda}_1 \le CK$ at the origin which will give us the required estimate (\ref{mainbd}).  Hence we may and do assume that $\tilde{\lambda}_1\gg K$ at this point.

We now differentiate $\tilde{H}$ at the origin, and as before, we use
 subscripts $k$ and $\ov{\ell}$ to denote the partial derivatives
 $\partial/\partial z^k$ and $\partial/\partial \ov{z}^{\ell}$. We have
\begin{equation}
\begin{aligned} \label{eq:Hk} \ti{H}_k &= \frac{\ti{\lambda}_{1,k}}{\ti{\lambda}_1} + \phi'(\alpha^{p\bar q}u_pu_{\bar qk} +
  \alpha^{p\bar q}u_{pk}u_{\bar q} + (\alpha^{p\bar q})_ku_pu_{\bar
q})  +
\psi' u_k \\
&= \frac{\ti{\lambda}_{1,k}}{\tilde{\lambda}_1} + \phi'(u_pu_{\bar pk} +
 u_{pk}u_{\bar p} + (\alpha^{p\bar q})_ku_pu_{\bar
q})  +
\psi' u_k \\
&= \frac{\ti{\lambda}_{1,k}}{\tilde{\lambda}_1} + \phi'V_k  +
\psi' u_k, \quad \textrm{for } V_k : = u_pu_{\bar pk} +
 u_{pk}u_{\bar p} + (\alpha^{p\bar q})_ku_pu_{\bar
q}.\end{aligned}\end{equation}

Differentiating once more,
\[ \begin{aligned}\ti{H}_{k\bar k} &= \frac{\ti{\lambda}_{1,k\bar k}}{\lambda_1} -
\frac{|\ti{\lambda}_{1,k}|^2}{\lambda_1^2}
+ \phi'\Big[u_pu_{\bar pk\bar k} + u_{\bar p}u_{pk\bar k} +
\sum_p|u_{pk}|^2 + \sum_p|u_{\bar pk}|^2 \\
& + (\alpha^{p\bar q})_{\bar k}u_pu_{\bar
    qk} + (\alpha^{p\bar q})_{\bar k}u_{pk}u_{\bar q} + (\alpha^{p\bar
    q})_{k\bar k}u_pu_{\bar q} + (\alpha^{p\bar q})_k(u_pu_{\bar q\bar
    k} + u_{p\bar k}u_{\bar q})\Big] \\
  &\quad+ \phi''|V_k|^2 + \psi''|u_k|^2 + \psi' u_{k\bar k},
  \end{aligned}\]
where we use the convention that we sum in all repeated indices except the free index $k$.

Since $(4K)^{-1}<\phi' < (2K)^{-1}$, we can absorb all the terms involving
$\alpha$ into the squared terms up to a constant, i.e. we have
\begin{equation}
\begin{aligned} \label{eq:b6} \ti{H}_{k\bar k} &\geq \frac{\ti{\lambda}_{1,k\bar k}}{\lambda_1} -
\frac{|\ti{\lambda}_{1,k}|^2}{\lambda_1^2}
+ \phi'(u_pu_{\bar pk\bar k} + u_{\bar p}u_{pk\bar k}) +
\frac{1}{5K} \sum_p(|u_{pk}|^2 + |u_{\bar pk}|^2) \\
  &\quad+ \phi''|V_k|^2 + \psi''|u_k|^2 + \psi' u_{k\bar k} - C.
 \end{aligned} \end{equation}
The constant $C$ denotes a constant that may change from line to line,
but it does not depend on the parameters $D_1,D_2$ that we are yet to choose.

\subsection*{$\bullet$ Calculation of $\widetilde{\lambda}_{1,k\bar k}$}
Let us now compute the derivatives of $\ti{\lambda}_1$. We have the
following general formulas for the derivatives of the eigenvalue
$\lambda_i$ of complex $n\times n$ matrices at a diagonal
matrix with distinct real eigenvalues
(see for instance Spruck~\cite{Spruck} in the case of matrices with real entries):
\begin{equation}
\begin{split}
\label{deriv} \lambda_i^{pq} = {} & \delta_{pi}\delta_{qi} \\
 \lambda_i^{pq,rs} = {} &
(1-\delta_{ip})\frac{\delta_{iq}\delta_{ir}\delta_{ps}}{\lambda_i -
  \lambda_p} +
(1-\delta_{ir})\frac{\delta_{is}\delta_{ip}\delta_{rq}}{\lambda_i -
  \lambda_r},
  \end{split}
  \end{equation}
where $\lambda_i^{pq}$ denotes the derivative with respect to the
$(p,q)$-entry $A_q^p$ of the matrix $A$, as a complex variable.

Denoting by $\ti{\lambda}_1$ the largest eigenvalue of the endomorphism $\ti{A}$
again, we have, using \eqref{deriv},
\begin{equation}
 \begin{aligned}\label{eq:b10} \ti{\lambda}_{1,k} &=
   \ti{\lambda}_1^{pq}\nabla_k(\ti{A}^p_q) \\
&= \nabla_k (\ti{A}^1_1) \\
&= \nabla_k g_{1\bar 1} - \nabla_k B^1_1 \\
&= g_{1\bar 1k} + (\alpha^{1\bar 1})_kg_{1\bar 1},
\end{aligned}\end{equation}
since $\nabla_k B^1_1 = 0$ at the origin.
Here we computed using covariant derivatives with respect to the Chern
connection of $\alpha$, which
makes the positivity of certain terms more apparent when we take second derivatives:
\begin{equation}
 \begin{aligned}\label{eq:z1} \ti{\lambda}_{1,k\bar k} &= \ti{\lambda}_1^{pq}\nabla_{\bar k}\nabla_k  \ti{A}^p_q +
\ti{\lambda}_1^{pq,rs} (\nabla_k \ti{A}^p_q) (\nabla_{\bar k} \ti{A}^r_s) \\
& = \nabla_{\ov{k}} \nabla_k g_{1\ov{1}}  +
\ti{\lambda}_1^{pq,rs} (\nabla_k \ti{A}^p_q) (\nabla_{\bar k} \ti{A}^r_s),
 \end{aligned}\end{equation}
where we used \eqref{deriv} and the fact that $\nabla_{\ov{k}}\nabla_k B^1_1=0$ at the origin.
To rewrite this in terms of partial derivatives, note first that
\begin{equation} \label{eq:fdg}
\begin{split}
\nabla_k g_{1\ov{1}} = {} & g_{1\ov{1}k} - \Gamma_{k1}^m g_{m\ov{1}} \\
 \nabla_{\bar k}\nabla_k g_{1\bar 1} =  {} &   g_{1\bar 1k\bar k} - (\partial_{\ov{k}} \Gamma_{k1}^m) g_{m\ov{1}} -  \Gamma_{k1}^m g_{m\ov{1}\ov{k}} -
 \ov{\Gamma^q_{k 1}} g_{1\ov{q} k} +  \ov{\Gamma^q_{k 1}}  \Gamma^m_{k1} g_{m\ov{q}} \\
 = {} & g_{1\bar 1k\bar k} + O( \sum_m
|g_{1\bar mk}|  + \lambda_1). \\
\end{split}
\end{equation}
In addition we have
\begin{equation}
\begin{aligned} \label{eq:l1}
\ti{\lambda}_1^{pq,rs}(\nabla_k \ti{A}^p_q) (\nabla_{\bar k} \ti{A}^r_s) &= \sum_{p > 1}
\frac{(\nabla_k \ti{A}^p_1)(\nabla_{\bar k} \ti{A}^1_p) + (\nabla_k \ti{A}^1_p)
  (\nabla_{\bar k} \ti{A}^p_1)}{\lambda_1 -
  \ti{\lambda}_p}\\
&=  \sum_{p > 1}
\frac{(\nabla_k g_{1\ov{p}}+\Gamma^q_{k1}B^p_q)(\nabla_{\bar k} g_{p\ov{1}}) + (\nabla_k g_{p\ov{1}}-\Gamma^1_{kq}B^q_p)
  (\nabla_{\bar k} g_{1\ov{p}})}{\lambda_1 -
  \ti{\lambda}_p}\\
&\geq \frac{1}{2}\sum_{p > 1} \frac{ |\nabla_k g_{1\bar p}|^2 + |\nabla_k g_{p\bar
    1}|^2}{\lambda_1 - \ti{\lambda}_p}-C,\end{aligned}\end{equation}
where we used \eqref{bdd} and \eqref{deriv}.
Recall that thanks to our choice of $B$ we have $\sum_i \ti{\lambda}_i>-1$, which implies $(\lambda_1 -
\ti{\lambda}_p)^{-1}\geq (n\lambda_1+1)^{-1}$ for $p > 1$, and so
\[ \ti{\lambda}_1^{pq,rs}(\nabla_k \ti{A}^p_q) (\nabla_{\bar k} \ti{A}^r_s) \geq
\frac{1}{2(n\lambda_1+1)} \sum_{p > 1} \Big( |\nabla_k g_{1\bar p}|^2 + |\nabla_k g_{p\bar
    1}|^2\Big)-C. \]
To rewrite this in terms of partial derivatives, note that
\[ \nabla_k g_{1\bar p} = g_{1\bar pk} - \Gamma^\ell_{k1}g_{\ell\bar p} =
g_{1\bar pk} + O(\lambda_1), \]
where we made use of the fact that $\sum_i\lambda_i>0$ to conclude that $|\lambda_i|\leq (n-1)\lambda_1$ for all $i$.
It follows, since we assume $\lambda_1 > 1$, that
\[ \ti{\lambda}_1^{pq,rs}(\nabla_k \ti{A}^p_q) (\nabla_{\bar k} \ti{A}^r_s) \geq
\frac{1}{4n\lambda_1} \sum_{p > 1} \Big( |g_{1\bar pk}|^2 + |g_{p\bar
  1k}|^2\Big) - C\lambda_1. \]

Combining this with \eqref{eq:z1} and \eqref{eq:fdg} we obtain
\begin{equation}
\begin{aligned} \label{eq:b1}\ti{\lambda}_{1,k\bar k} &\geq g_{1\bar 1k\bar k} +
\frac{1}{4n\lambda_1} \sum_{p > 1} \Big( |g_{1\bar pk}|^2 + |g_{p\bar
  1k}|^2\Big) - C( \sum_m |g_{1\bar mk}|  + \lambda_1) \\
&\geq g_{1\bar 1k\bar k} + \frac{1}{8n\lambda_1}
\sum_{p > 1} \Big( |g_{1\bar pk}|^2 + |g_{p\bar
  1k}|^2\Big) - C(|g_{1\bar 1k}| + \lambda_1).\end{aligned}\end{equation}

Rewriting $g$ in terms of $u$, we have
\begin{equation}
\begin{aligned} \label{eq:b2} g_{1\bar 1k\bar k} &= \chi_{1\bar 1k\bar k} + u_{1\bar 1k\bar k} +
W_{1\bar 1k\bar k} \\
&= \chi_{1\bar 1k\bar k} + u_{k\bar k1\bar 1} + W_{1\bar 1k\bar k} \\
&= \chi_{1\bar 1k\bar k} -
\chi_{k\bar k1\bar 1} + g_{k\bar k1\bar 1} - W_{k\bar k1\bar 1}
+W_{1\bar 1k\bar k},
\end{aligned}\end{equation}
and so
\begin{equation}
\begin{aligned} \label{eq:j2} F^{kk}\widetilde{\lambda}_{1,k\bar k} &\geq F^{kk}g_{k\bar k1\bar 1}
+ F^{kk}(W_{1\bar 1k\bar k} - W_{k\bar k1\bar 1}) \\
&\quad + \frac{1}{8n\lambda_1}
\sum_{p > 1} F^{kk}\Big( |g_{1\bar pk}|^2 + |g_{p\bar
  1k}|^2\Big) \\
&\quad - C(F^{kk}|g_{1\bar 1k}| + \lambda_1\mathcal{F}),\end{aligned}\end{equation}
where $F^{pq}$ denotes the partial derivative of the function $F(A)$ with respect to the $(p,q)$-entry of the matrix $A$ (as explained earlier), and we have set $\mathcal{F} = \sum_k F^{kk}$.
Observe that, thanks to \eqref{use} and \eqref{deriv}, at the origin we have that $F^{pq}$ vanishes whenever $p\neq q$, while on the other hand $F^{kk}=f_k$, using the notation from section \ref{sec:bg}.
Recall from the last assertion of Proposition \ref{prop:subsol} that
\begin{equation}
\label{eq:mcFlb}
 \mathcal{F} \ge \kappa>0,\end{equation}
for a uniform $\kappa>0$.

\subsection*{$\bullet$ The term $F^{kk}g_{k\bar k1\bar 1}$}
We now differentiate the equation $F(A)=h$, using covariant
derivatives to simplify a term
that appears below. Applying $\nabla_i$, we obtain
\[ F^{pq} \nabla_i g_{q\ov{p}} = h_i, \]
namely,
\begin{equation}
 \label{eq:b3} F^{kk} g_{k\bar ki} + F^{kk}(\alpha^{k\bar
  k})_ig_{k\bar k} = h_i. \end{equation}

Applying $\nabla_{\ov{i}}$ and setting $i=1$,
\begin{equation}
 \label{eq:F11}
F^{pq,rs} \nabla_1g_{q\bar p}\nabla_{\bar 1}g_{s\bar r} + F^{kk}
\nabla_{\bar 1}\nabla_1 g_{ k\bar k} = h_{1\bar 1}.
\end{equation}
To rewrite this using partial derivatives, note that
\begin{equation}
\begin{split}
\nabla_{\bar 1}\nabla_1 g_{k\bar k} = {} &  g_{k\bar k1\bar 1} - (\partial_{\ov{1}} \Gamma_{1k}^m) g_{m\ov{k}} -  \Gamma_{1k}^m g_{m\ov{k}\ov{1}} -
 \ov{\Gamma^q_{1 k}} g_{k\ov{q} 1} + \ov{\Gamma^q_{1k}}  \Gamma^m_{1k} g_{m\ov{q}} \\
={} &
 g_{k\bar k1\bar 1} -  2\textrm{Re} \left( \ov{\Gamma^q_{1k}} g_{k\ov{q} 1} \right) +
O( \lambda_1).
\end{split}
\end{equation}
By rewriting $g$ in terms of $u$, we have
\[ g_{k\ov{q}1} = g_{1\ov{q}k} + \chi_{k\ov{q}1} - \chi_{1\ov{q}k} + W_{k\ov{q}1} - W_{1\ov{q}k}, \]
and hence
\[ \nabla_{\bar 1}\nabla_1 g_{k\bar k} = g_{k\bar k1\bar 1} -2\textrm{Re} \left( \ov{\Gamma^q_{1k}} (W_{k\ov{q}1} - W_{1\ov{q}k}) \right) +
O(\sum_q |g_{1\bar q k}|  + \lambda_1). \]

Returning to (\ref{eq:F11}), and making use of  (\ref{eq:mcFlb}), we obtain
\begin{equation}\label{eq:j1}
\begin{split}
F^{kk} g_{k\bar k1\bar 1} \geq {} &  -F^{pq,rs}\nabla_1g_{q\bar
  p}\nabla_{\bar 1}g_{s\bar r} - CF^{kk} \sum_q |g_{1\bar qk}|  - C\mathcal{F}\lambda_1 \\
 & {} + 2\textrm{Re} \left( F^{kk} \ov{\Gamma^q_{1k}} (W_{k\ov{q}1} -
   W_{1\ov{q}k}) \right).
 \end{split}
  \end{equation}
To bound the term involving $W_{1\bar qk}$, we observe that
\begin{equation}
 \label{eq:W1qk}
|F^{kk}W_{1\bar qk}| \leq  C(\mathcal{F}\lambda_1 +
\sum_pF^{kk}|u_{pk}|). \end{equation}
To see (\ref{eq:W1qk}), we use that $\lambda_1>K$ to bound the terms involving the gradient of $u$ that arise from $W_{1\ov{q}k}$ when taking the $\de/\de z^k$ derivative of $W_{1\ov{q}}=(\tr{\alpha}{Z})\alpha_{1\ov{q}}-(n-1)Z_{1\ov{q}}$. For the  term involving $W_{k \bar q1}$  note that
\[ \begin{aligned}\nabla_{\bar 1}\nabla_1 W_{k\bar k} = {}&  \nabla_{\bar 1}(W_{k\bar k1}
- \Gamma^q_{1k}W_{q\bar k}) \\
= {} & W_{k\bar k1\bar 1} - \overline{\Gamma^q_{1k}}W_{k\bar q1} -
\Gamma^q_{1k}W_{q\bar k\bar 1} \\
{} & - (\Gamma^q_{1k})_{\bar 1} W_{q\bar k}
+ \overline{\Gamma^p_{1k}}\Gamma^q_{1k} W_{q\bar p}. \end{aligned}\]
In particular
\[ 2\mathrm{Re}(\overline{\Gamma^q_{1k}}W_{k\bar q1}) =  W_{k\bar k1\bar 1} - \nabla_{\bar
  1}\nabla_1 W_{k\bar k}+ O(K^{1/2}). \]
Using this (and that we can assume $\lambda_1 > K$), we have
\begin{equation}
\begin{aligned} \label{eq:Wkq1}
2F^{kk}\mathrm{Re}(\overline{\Gamma^q_{1k}}W_{k\bar q1})  &= F^{kk}
W_{k\bar k1\bar 1} -F^{kk} \nabla_{\bar  1}\nabla_1 W_{k\bar k}   +
O(\mathcal{F}\lambda_1).\end{aligned}\end{equation}
Combining \eqref{eq:j1}, \eqref{eq:W1qk} and \eqref{eq:Wkq1} gives
\begin{equation}
 \begin{aligned}\label{eq:penu}
F^{kk}g_{k\bar k1\bar 1} &\geq - F^{pq,rs}\nabla_1g_{q\bar p}
\nabla_{\bar 1}g_{s\bar r} + F^{kk}W_{k\bar k1\bar 1} -
F^{kk}\nabla_{\bar 1}\nabla_1 W_{k\bar k} \\
&\quad - C\left( F^{kk}\sum_q |g_{1\bar qk}| + F^{kk}\sum_p |u_{pk}| +
  \mathcal{F}\lambda_1\right). \end{aligned}\end{equation}
Going back to \eqref{eq:j2}, using the square terms there to control the terms in (\ref{eq:penu}) involving $|g_{1\ov{q}k}|$ for $q\neq 1$, we obtain
\begin{equation} \label{eq:b5}
\begin{split}
F^{kk}\ti{\lambda}_{1,k\bar k} \geq {} & -F^{pq,rs}\nabla_1g_{q\bar
  p}\nabla_{\bar 1}g_{s\bar r} + F^{kk}(W_{1\bar 1k\bar k} -
\nabla_{\bar 1}\nabla_1 W_{k\bar k})
\\
&\quad - C(F^{kk}|g_{1\bar 1k}| + \sum_pF^{kk}|u_{pk}| + \lambda_1\mathcal{F}).
   \end{split}
 \end{equation}

\subsection*{$\bullet$ The term $F^{kk}\nabla_{\bar 1}\nabla_1 W_{k\bar k}$}
Using \eqref{eq:cz} and \eqref{fk} we have
\begin{equation}
\begin{aligned} \label{eq:d1} F^{kk}\nabla_{\bar 1}\nabla_1 W_{k\bar k} &=
\frac{1}{n-1}\sum_k \sum_{i\ne k} \widetilde{F}^{ii}
\nabla_{\bar 1}\nabla_1 W_{k\bar k} \\
&= \frac{1}{n-1}\sum_i \widetilde{F}^{ii}\sum_{k\ne i}
\nabla_{\bar 1}\nabla_1 W_{k\bar k} \\
&=  \widetilde{F}^{ii} \nabla_{\bar 1}\nabla_1 Z_{i\bar i},\end{aligned}\end{equation}
Recall from \eqref{eq:fk1} and \eqref{eq:fk2} that $\tilde{F}^{11} = \tilde{f}_1$ is ``large'', equivalent to $F^{kk}=f_k$ for any $k>1$, while $\tilde{F}^{ii} =\tilde{f}_i$ for $i>1$ is ``small'', bounded by $F^{11}=f_1$.
We also recall that, as explained earlier, the crucial assumption on $Z_{i\ov{j}}$ implies that $\nabla_{\bar 1}\nabla_1 Z_{1\bar 1}$
does not contain the terms $u_{11\bar 1}, u_{11}$ or their complex
conjugates.

Hence, using the fact that $\sup_{i,j} |u_{i\ov{j}}| \le C \lambda_1$ and $\lambda_1 \ge K$,
\begin{equation}
 \begin{aligned}\label{eq:j5}
\lefteqn{
F^{kk}\nabla_{\bar 1}\nabla_1 W_{k\bar k} } \\
\leq {} & C\Big( \widetilde{F}^{11} \sum_{k>1}( |u_{k1}|+|u_{k\ov{1}1}|) +  F^{11} \sum_k (|u_{k1}|+|u_{k\ov{1}1}|) + \lambda_1 \mathcal{F} \Big) \\
\leq {} & C\Big(F^{11}(|u_{11}| + |u_{1\bar 11}|) +  \sum_{k > 1}
F^{kk}(|u_{k1}| + |u_{k\bar 11}|) + \lambda_1\mathcal{F} \Big) \\
\leq {} & C\Big( \sum_p F^{kk}|u_{pk}| + F^{kk}|u_{1\bar 1k}| +
\lambda_1\mathcal{F}\Big). \end{aligned}\end{equation}
We also have
\[ u_{1\bar 1k} = g_{1\bar 1k} - \chi_{1\bar 1k} - W_{1\bar 1k}, \]
and so
\[ F^{kk}|u_{1\bar 1k}| \leq F^{kk}|g_{1\bar 1k}| + C\Big[
\sum_pF^{kk}|u_{pk}| + \lambda_1\mathcal{F}\Big]. \]
From \eqref{eq:j5} we then obtain
\begin{equation}
\begin{aligned} \label{eq:j6} F^{kk}\nabla_{\bar 1}\nabla_1 W_{k\bar k}
&\leq C\Big( \sum_p F^{kk}|u_{pk}| + F^{kk}|g_{1\bar 1k}| +
\lambda_1\mathcal{F}\Big). \end{aligned}\end{equation}

\subsection*{ $\bullet$ The term $F^{kk}W_{1\bar 1k\bar k}$}

Let us write $W_{1\bar 1} =
W^pu_p + W^{\bar p}u_{\bar p}$. We have
\[ (W^pu_p)_{k\bar k} = (W^p)_{k\bar k}u_p + (W^p)_k u_{p\bar k} +
(W^p)_{\bar k} u_{pk} + W^pu_{pk\bar k}, \]
and so since we can assume that $\lambda_1 \gg K$, we have
\begin{equation}
 \label{eq:b4} F^{kk}(W^pu_p)_{k\bar k} \geq W^pF^{kk}u_{k\bar kp} -C\Big[
\sum_pF^{kk}|u_{pk}| + \lambda_1\mathcal{F}\Big]. \end{equation}
Using \eqref{eq:b3} we have
\[ \begin{aligned}|F^{kk}u_{k\bar kp}| &= |F^{kk}(g_{k\bar kp} - \chi_{k\bar kp} -
W_{k\bar kp})| \\
&\leq CF^{kk}|g_{k\bar k}| + C\mathcal{F}  + |F^{kk}W_{k\bar kp}| \\
&\leq CF^{kk}|u_{k\bar k}| +  CK^{1/2}\mathcal{F} + |F^{kk}W_{k\bar kp}|.\end{aligned}\]
To deal with this last term, note that thanks to \eqref{eq:cz}, as in
\eqref{eq:d1}
\[ \begin{split}
F^{kk}W_{k\bar kp} &= F^{kk} \nabla_p W_{k\bar k} +
O(K^{1/2}\mathcal{F}) \\
&= \widetilde{F}^{ii}\nabla_p Z_{i\bar i} + O(K^{1/2}\mathcal{F}),
\end{split} \]
and using the crucial assumption on $Z_{i\ov{j}}$, as explained earlier,
we see that $\nabla_p Z_{1\ov{1}}$ is independent of $u_{11},
u_{1\ov{1}}$. It follows that these Hessian terms
 can appear only with the ``small'' coefficients $\widetilde{F}^{ii}$
with $i > 1$. We obtain
\[ \begin{aligned}|F^{kk}W_{k\bar kp}| \leq {} & C \Big( \tilde{F}^{11} \sum_{k>1}(|u_{kp}| + |u_{\ov{k}p}|) + F^{11} \sum_k (|u_{kp}| + |u_{\ov{k}p}|) + K^{1/2} \mathcal{F} \Big) \\
\leq {}  & C\big( \sum_pF^{kk}|u_{pk}| +
\sum_pF^{kk}|u_{p\bar k}| + K^{1/2}\mathcal{F}   \big), \end{aligned}\]
and so
\begin{equation}
\begin{aligned}  \label{eq:extra} |F^{kk}u_{k\bar kp}| &\leq C\Big[ \sum_pF^{kk}|u_{pk}| + \sum_pF^{kk}|u_{p\bar k}| +
K^{1/2}\mathcal{F}\Big]. \end{aligned}\end{equation}
From \eqref{eq:b4} we then have (using $\lambda_1 \gg K$):
\[ F^{kk}(W^pu_p)_{k\bar k} \geq - C\Big[ \sum_pF^{kk}|u_{pk}| +
\sum_pF^{kk}|u_{p\bar k}| + \lambda_1\mathcal{F} \Big]. \]
A similar argument gives the same estimate for $F^{kk}(W^{\bar
  p}u_{\bar p})_{k\bar k}$, and this completes the required estimate
for $F^{kk} W_{1\ov{1} k\ov{k}}$:
\[ F^{kk}W_{1\bar 1 k\bar k} \geqslant -C\Big[ \sum_p F^{kk}|u_{pk}| +
\lambda_1\mathcal{F}\Big]. \]

Putting together this last inequality and \eqref{eq:j6} into \eqref{eq:b5} we obtain
\[ F^{kk}\ti{\lambda}_{1,k\bar k} \geq -F^{pq,rs}\nabla_1g_{q\bar
  p}\nabla_{\bar 1}g_{s\bar r} - C\Big[ F^{kk}|g_{1\bar 1k}| +
\sum_pF^{kk}|u_{pk}| + \lambda_1\mathcal{F}\Big]. \]
We now use this in Equation~\eqref{eq:b6}, to give
\[ \begin{aligned}F^{kk}\ti{H}_{k\bar k} &\geq \frac{-F^{pq,rs}\nabla_1g_{q\bar
    p}\nabla_{\bar 1}g_{s\bar r}}{\lambda_1} - \frac{F^{kk}|\ti{\lambda}_{1,k}|^2}{\lambda_1^2} + F^{kk}\phi'(u_pu_{\bar pk\bar k} +
u_{\bar p}u_{pk\bar k}) \\
&\quad + \sum_p\frac{F^{kk}}{5K}(|u_{pk}|^2 + |u_{\bar pk}|^2) +
\phi''F^{kk}|V_k|^2 + \psi''F^{kk}|u_k|^2 + \psi'F^{kk}u_{k\bar k} \\
&\quad - C\Big[ F^{kk}\lambda_1^{-1} |g_{1\bar 1k}| +
\sum_pF^{kk}\lambda_1^{-1} |u_{pk}| + \mathcal{F}\Big]. \end{aligned}\]

We can use \eqref{eq:extra} and the fact that $\phi' < (2K)^{-1}$ to
bound the terms involving $u_{\bar pk\bar k}, u_{pk\bar k}$:
\[ \sum_p\phi' F^{kk} |u_pu_{\bar pk\bar k}| \leq C\Big[ \frac{1}{K^{1/2}} \sum_pF^{kk}|u_{pk}| +
\frac{1}{K^{1/2}} \sum_pF^{kk}|u_{p\bar k}| + \mathcal{F}\Big], \]
which in turn can be controlled by the good squared terms $|u_{pk}|^2
+ |u_{p\bar k}|^2$ at the cost of an extra multiple of $\mathcal{F}$. In
addition, since we assume $\lambda_1 \gg K$, we can control the
$F^{kk}\lambda_1^{-1}|u_{pk}|$ term  in the same way. We therefore have
\begin{equation}
 \begin{aligned}\label{eq:b7} 0 \ge F^{kk}\ti{H}_{k\bar k} \geq {} &  \frac{-F^{pq,rs}\nabla_1g_{q\bar
    p}\nabla_{\bar 1}g_{s\bar r}}{\lambda_1} -
\frac{F^{kk}|\ti{\lambda}_{1,k}|^2}{\lambda_1^2} \\
& {} + \sum_p\frac{F^{kk}}{6K}(|u_{pk}|^2
  + |u_{p\bar k}|^2) +
\phi'' F^{kk} |V_k|^2 \\ &  {} + \psi''F^{kk} |u_k|^2 + \psi'F^{kk} u_{k\bar k}
 - C\Big[
F^{kk}\lambda_1^{-1}|g_{1\bar 1k}| + \mathcal{F}\Big]. \end{aligned}\end{equation}

We now deal with two cases separately, as was done in
Hou-Ma-Wu~\cite{HMW}, depending on a small constant $\delta =
\delta_{D_1,D_2}>0$ to be determined shortly, and which will depend on the constants $D_1$ and $D_2$.

\bigskip\noindent{\bf Case 1.} Assume $\delta\lambda_1 \geq
-\lambda_n$.  Define the set
\[ I = \{ i\,:\, F^{ii} > \delta^{-1}F^{11}\,\}. \]
From (\ref{eq:Hk}) and the fact that  $\ti{H}_k=0$ at the maximum, we get
\begin{equation}
\begin{aligned}\label{eq:b8} -\sum_{k\not\in I} \frac{F^{kk}|\ti{\lambda}_{1,k}|^2}{\lambda_1^2} &=
-\sum_{k\not\in I} F^{kk}\left|\phi' V_k + \psi'u_k\right|^2 \\
&\geq -2\phi'^2 \sum_{k\not\in I} F^{kk}|V_k|^2 -
2\psi'^2\sum_{k\not\in I} F^{kk}|u_k|^2 \\
&\geq -\phi'' \sum_{k\not\in I} F^{kk}|V_k|^2 - 2\psi'^2\delta^{-1}F^{11}K. \end{aligned}\end{equation}

For $k\in I$ we have in the same way
\begin{equation}
 \label{eq:b9} -2\delta \sum_{k\in I} \frac{F^{kk}|\ti{\lambda}_{1,k}|^2}{\lambda_1^2} \geq
- 2\delta \phi'' \sum_{k\in I} F^{kk}|V_k|^2 - 4\delta \psi'^2\sum_{k\in I}
F^{kk}|u_k|^2. \end{equation}
We wish to use some of the good $\psi''F^{kk}|u_k|^2$ term in \eqref{eq:b7} to control the last term in \eqref{eq:b9}. For
this we assume that $\delta$ is chosen so small (depending on $\psi$,
i.e. on $D_1,D_2$ and the maximum of $|u|$), such that
\begin{equation} \label{eq:definedelta}
4\delta\psi'^2 < \frac{1}{2}\psi''
\end{equation}
Since $\psi''$ is strictly positive, such a $\delta > 0$ exists.

Using this together with \eqref{eq:b8}, \eqref{eq:b9} in \eqref{eq:b7}, we have
\begin{equation}
 \begin{aligned}\label{eq:b16} 0 &\geq \frac{-F^{pq,rs}\nabla_1g_{q\bar
    p}\nabla_{\bar 1}g_{s\bar r}}{\lambda_1} -
(1-2\delta) \sum_{k\in I} \frac{F^{kk}|\ti{\lambda}_{1,k}|^2}{\lambda_1^2}
\\
\quad& + \sum_p\frac{F^{kk}}{6K}(|u_{pk}|^2  + |u_{p\bar k}|^2) +
\frac{1}{2} \psi'' F^{kk}|u_k|^2 + \psi'F^{kk}u_{k\bar k}\\
&\quad - 2\psi'^2\delta^{-1}F^{11}K - C\Big[
F^{kk}\lambda_1^{-1}|g_{1\bar 1k}| + \mathcal{F}\Big]. \end{aligned}\end{equation}

To deal with the first two terms, note that (as in \cite[Equation
(67)]{Sz})
the concavity of the
operator $F$ implies
\begin{equation}
\label{eq:b15} - F^{pq,rs}\nabla_1g_{q\bar p}\nabla_{\bar 1}g_{s\bar r} \geq
 \sum_{k\in I} \frac{F^{kk}- F^{11}}{\lambda_1 - \lambda_k}
|\nabla_1g_{k\bar 1}|^2, \end{equation}
where note that the denominator involves $\lambda_k$ instead of $\ti{\lambda}_k$ because we are evaluating $F$ at $A$. We also remark that
the denominator on the right hand side does not vanish, because the assumption $k\in I$ implies that $F^{kk}>F^{11}$, which implies that $\lambda_k<\lambda_1$ because $f$ is symmetric.
By definition, for $k\in I$ we have $F^{11}\leq \delta F^{kk}$, and the assumption that $\delta \lambda_1 \geq -\lambda_n$ implies
\[ \frac{1- \delta}{\lambda_1 - \lambda_k} \geq
\frac{1-2\delta}{\lambda_1}. \]
It follows that
\begin{equation}
 \begin{aligned}\label{eq:b14} \sum_{k\in I} \frac{F^{kk}- F^{11}}{\lambda_1 - \lambda_k}
|\nabla_1g_{k\bar 1}|^2 &\geq \sum_{k\in I} \frac{(1-\delta
 )F^{kk}}{\lambda_1 - \lambda_k} |\nabla_1g_{k\bar 1}|^2 \\
&\geq \frac{1-2\delta}{\lambda_1} \sum_{k\in I} F^{kk}|\nabla_1g_{k\bar
  1}|^2.\end{aligned}\end{equation}
Combining this with \eqref{eq:b16} and \eqref{eq:b15}, we then obtain
\begin{equation}
\begin{aligned} \label{eq:j9}
0 &\geq
(1-2\delta) \sum_{k\in I} \frac{F^{kk}(|\nabla_1g_{k\bar 1}|^2 -
  |\ti{\lambda}_{1,k}|^2)}{\lambda_1^2}
\\
\quad& + \sum_p\frac{F^{kk}}{6K}(|u_{pk}|^2  + |u_{p\bar k}|^2) +
\frac{1}{2} \psi'' F^{kk}|u_k|^2 + \psi'F^{kk}u_{k\bar k}\\
&\quad - 2\psi'^2\delta^{-1}F^{11}K - C\Big[
F^{kk}\lambda_1^{-1}|g_{1\bar 1k}| + \mathcal{F}\Big]. \end{aligned}\end{equation}

We wish to obtain a lower bound for the first term in \eqref{eq:j9}.  We make the following claim.

\bigskip
\noindent
{\bf Claim.}  \ For any $\ve>0$ there exists a constant $C_{\ve}$ such that
\begin{equation}
 \begin{aligned}\label{eq:squares}
\sum_{k\in I}\frac{F^{kk}|\nabla_1 g_{k\bar 1}|^2}{\lambda_1^2}
\geq {} & \sum_{k\in I} \frac{F^{kk}|\ti{\lambda}_{1,k}|^2}{\lambda_1^2} -
\sum_p\frac{F^{kk}}{12K} (|u_{p\bar k}|^2 + |u_{pk}|^2) \\ {} & + C_\epsilon \psi'
F^{kk} |u_k|^2
+ \epsilon C\psi' \mathcal{F} - C\mathcal{F},
\end{aligned}\end{equation}
as long as $\lambda_1 /K$ is sufficiently large compared to $\psi'$ (the constants $D_1$, $D_2$ of $\psi$ will  be chosen uniformly later).

\begin{proof}[Proof of Claim]
First, we compare $\nabla_1g_{k\bar 1}$ to $\widetilde{\lambda}_{1,k}$.  We have
\[ \begin{aligned}\nabla_1 g_{k\bar 1} &= g_{k\bar 11} - \Gamma^p_{1k}g_{p\bar 1} \\
&= \chi_{k\bar 11} + u_{k\bar 11} + W_{k\bar 11} + O(\lambda_1) \\
&= \chi_{k\bar 11} + u_{1\bar 1k} + W_{k\bar 11} + O(\lambda_1) \\
&= g_{1\bar 1k} - W_{1\bar 1k} + W_{k\bar 11} + O(\lambda_1)\\
&= \ti{\lambda}_{1,k} - W_{1\bar 1k} + W_{k\bar 11} + O(\lambda_1), \end{aligned}\]
absorbing bounded terms into $O(\lambda_1)$ and using
\eqref{eq:b10}. It follows that for any $k$, without summing,
\begin{equation}
 \begin{aligned}\label{eq:b11} |\nabla_1g_{k\bar 1}|^2 \geq {} &
   |\ti{\lambda}_{1,k}|^2 -
C\Big[
 |\ti{\lambda}_{1,k}|(|W_{1\bar 1k}| + |W_{k\bar
  11}|) \\
& + \lambda_1 |\ti{\lambda}_{1,k}| + |W_{1\bar 1k}|^2 + |W_{k\bar 11}|^2 +
\lambda_1^2\Big].\end{aligned}\end{equation}

\subsection*{$\bullet$ The terms in \eqref{eq:b11} involving $W$}
Note that if $k\in I$ then $k\ne 1$, and so from \eqref{eq:cz} we have
\[ W_{k\bar 11} = (\mathrm{tr}_\alpha Z)\alpha_{k\bar 11} -
(n-1)Z_{k\bar 11}. \]
Our basic assumption for $Z$ implies that $Z_{k\bar 11}=\nabla_1
Z_{k\bar 1} + O(Z)$ does not contain the Hessian terms $u_{11}$ or $u_{\bar 1\bar
  1}$. It follows that $W_{k\bar 11}$ and its complex conjugate do not
contain these Hessian terms.
The term $W_{1\bar 1k}$ and its complex conjugate also do
not contain the Hessian terms
$u_{11}$ or $u_{\bar 1\bar 1}$ since each Hessian term must contain a
$k$-derivative. To simplify the formulas, let us write
\[ U = \sum_{\substack{p > 1\\q\geqslant 1}} |u_{pq}|. \]
It follows that
\begin{equation}\label{eq:Wk11}
  |W_{k\bar 11}| + |W_{1\bar 1k}| \leqslant C\left( \lambda_1 +
    U\right),
\end{equation}
and so the terms $|W_{k\bar 11}|^2 + |W_{1\bar 1k}|^2$ in \eqref{eq:b11} can be bounded by $C( \lambda_1^2 +
  U^2).$

We now use these to estimate the negative terms in
\eqref{eq:b11}. Using that $\ti{H}_k=0$ together with \eqref{eq:Wk11}
 we have
\begin{equation}
 \begin{aligned}\label{eq:j8} &|\ti{\lambda}_{1,k}|(|W_{1\bar 1k}| + |W_{k\bar 11}|) \\
&=
\lambda_1 \Big|\phi'(u_ru_{\bar rk} + u_{rk}u_{\bar r} +
(\alpha^{r\bar s})_ku_r u_{\bar s})
+ \psi' u_k\Big| (|W_{1\bar 1k}| + |W_{k\bar 11}|)\\
&\leq \frac{C\lambda_1}{2K^{1/2}}( \sum_r|u_{\bar rk}| + \sum_r|u_{rk}| + K^{1/2})\left(
  \lambda_1 + U\right) \\
&\quad + C \lambda_1 |\psi'| |u_k| \left(\lambda_1 + U\right).
 \end{aligned}\end{equation}
We have
\begin{equation}
 \begin{aligned} \label{eq:tres}
\lefteqn{
\frac{C\lambda_1}{2K^{1/2}} (\sum_r |u_{\ov{r}k}| + \sum_r |u_{rk}| + K^{1/2}) \lambda_1 } \\
&\leq \frac{C\lambda_1^2}{2K^{1/2}} \sum_r |u_{\ov{r}k}| +
\frac{C\lambda_1^2}{2K^{1/2}}U + C\lambda_1^2 ,
  \end{aligned}\end{equation}
 and
  \begin{equation}
  \begin{aligned}\label{eq:cuatro}
    \lefteqn{\frac{C\lambda_1}{2K^{1/2}} \left(\sum_r |u_{\ov{r}k}| + \sum_r
      |u_{rk}| + K^{1/2}\right) U } \\ & \leq
  \frac{C\lambda_1^2}{K^{1/2}} U + \frac{C\lambda_1}{K^{1/2}} U^2 +
  C\lambda_1 U \\
&\leqslant \frac{C\lambda_1^2}{K^{1/2}} U + \frac{C\lambda_1}{K^{1/2}} U^2.
 \end{aligned}\end{equation}
Next, for any $\ve>0$, there exists a constant $C_{\ve}$ such that
 \begin{equation}
  \label{eq:uno}
 \begin{aligned}
 \lambda_1 |\psi'| |u_k| \lambda_1
 \leq {} & - \ve \lambda_1^2 \psi' - \lambda_1^2 C_{\ve} \psi'  |u_k|^2,
 \end{aligned}\end{equation}
 where we have used the fact that $\psi'<0$.  And,
 \begin{equation}
  \label{eq:dos} \lambda_1 | \psi'| |u_k| U \le -  \lambda_1 \psi'
  U^2 -  \lambda_1 \psi' |u_k|^2. \end{equation}
Combining \eqref{eq:j8} with  \eqref{eq:tres}, \eqref{eq:cuatro}, \eqref{eq:uno} and \eqref{eq:dos} we obtain
 \begin{equation}
  \begin{aligned}\label{eq:j8B} |\ti{\lambda}_{1,k}|(|W_{1\bar 1k}| + |W_{k\bar 11}|)
&\leq  C \bigg( \frac{\lambda_1^2}{2K^{1/2}}\sum_r |u_{\bar rk}| +
\frac{\lambda_1^2}{K^{1/2}}
  U + \frac{\lambda_1}{K^{1/2}}U^2 + \lambda_1^2 \\ &\quad
 - \lambda_1^2 C_\epsilon \psi' |u_k|^2  - \epsilon \lambda_1^2
\psi'- \lambda_1 \psi' |u_k|^2 - \lambda_1\psi' U^2 \bigg).
\end{aligned}\end{equation}
Using $\ti{H}_k = 0$ again,
\begin{equation}
\begin{aligned} \label{eq:b142}
  \frac{|\ti{\lambda}_{1,k}|}{\lambda_1} &=
|\phi'(u_pu_{\bar pk} + u_{pk}u_{\bar
  p} + (\alpha^{p\bar q})_ku_pu_{\bar q} ) + \psi' u_k|  \\
&\leq \frac{1}{2K^{1/2}} \sum_p |u_{\bar pk}| + \frac{1}{2K^{1/2}} U -
C_\epsilon \psi' |u_k|^2 - \epsilon \psi' + C. \end{aligned}\end{equation}
We then obtain
\[ \begin{aligned} \frac{|\nabla_1 g_{k\bar 1}|^2}{\lambda_1^2} &\geq
  \frac{|\ti{\lambda}_{1,k}|^2}{\lambda_1^2} - C\Big[
  \frac{1}{K^{1/2}}\sum_p |u_{p\bar k}| + \frac{1}{K^{1/2}}U \\
&\quad\quad + \frac{1}{\lambda_1 K^{1/2}}U^2 + 1 - \epsilon\psi' -
\frac{1}{\lambda_1}\psi'|u_k|^2 - \frac{1}{\lambda_1}\psi' U^2\Big] + C_\epsilon\psi'
|u_k|^2.
\end{aligned}\]
Summing over $k\in I$, we have
\begin{equation}\label{uno} \begin{aligned}
  \sum_{k\in I} \frac{F^{kk}|\nabla_1 g_{k\bar 1}|^2}{\lambda_1^2}  &
  \geq \sum_{k\in I}\frac{F^{kk}|\ti{\lambda}_{1,k}|^2}{\lambda_1^2} -
  C\Big[ \frac{1}{K^{1/2}}\sum_{\substack{p}} F^{kk}|u_{p\bar
      k}| + \frac{1}{K^{1/2}} \mathcal{F}U \\
&\quad + \frac{1}{\lambda_1 K^{1/2}} \mathcal{F}U^2 + \mathcal{F} - \epsilon
      \mathcal{F}\psi' - \frac{1}{\lambda_1}\psi' F^{kk}|u_k|^2 -
      \frac{1}{\lambda_1} \psi'\mathcal{F}U^2\Big]
      \\
&\quad +
      C_\epsilon\psi' F^{kk} |u_k|^2.
\end{aligned}\end{equation}
First, we use
\begin{equation}\label{due} \frac{C}{K^{1/2}}\sum_{\substack{p}} F^{kk}|u_{p\bar k}|\leq \frac{1}{12K} \sum_p F^{kk}|u_{p\ov{k}}|^2 +C\mathcal{F}.\end{equation}
Note that all $F^{kk}$ with $k > 1$ are comparable to
$\mathcal{F}$. It follows that
\begin{equation}\label{tre}  \frac{C}{K^{1/2}} \mathcal{F} U \leq \frac{1}{50K} \sum_p F^{kk}|u_{pk}|^2 + C
\mathcal{F}, \end{equation}
and
\begin{equation}\label{quattro}  \frac{C}{\lambda_1 K^{1/2}} \mathcal{F} U^2 \leq \frac{C}{\lambda_1
  K^{1/2}} \sum_p F^{kk}|u_{pk}|^2\leq \frac{1}{50K} \sum_p F^{kk}|u_{pk}|^2 . \end{equation}
As long as $\lambda_1 / K$ is sufficiently large
depending on $\psi'$ (i.e. depending on
$D_1,D_2$ which will be chosen later) we have
\begin{equation}\label{cinque} \frac{C}{\lambda_1} \psi'\mathcal{F}U^2\leq \frac{1}{50K} \sum_p F^{kk}|u_{pk}|^2.\end{equation}
and using \eqref{due}, \eqref{tre}, \eqref{quattro} and \eqref{cinque} in \eqref{uno} we finally obtain
\[\begin{aligned}
  \sum_{k\in I} \frac{F^{kk}|\nabla_1 g_{k\bar 1}|^2}{\lambda_1^2}  &
  \geq \sum_{k\in I}\frac{F^{kk}|\ti{\lambda}_{1,k}|^2}{\lambda_1^2} -
\sum_p\frac{F^{kk}}{12K} (|u_{p\bar k}|^2 + |u_{pk}|^2) \\
&\quad + C_\epsilon \psi'
F^{kk} |u_k|^2 + \epsilon C\psi' \mathcal{F} - C\mathcal{F}.
\end{aligned}\]
This completes the proof of the claim.
\end{proof}

 We now use the claim in
 \eqref{eq:j9} to
obtain
\begin{equation} \label{eq:star}
 \begin{aligned}0 &\geq \sum_p\frac{F^{kk}}{12K}(|u_{pk}|^2 + |u_{p\bar k}|^2) +
\frac{1}{2}\psi'' F^{kk}|u_k|^2 + \psi' F^{kk}u_{k\bar k} \\
&\quad - 2\psi'^2\delta^{-1}F^{11}K - C\big[
F^{kk}\lambda_1^{-1}|g_{1\bar 1k}| + \mathcal{F}\big] \\
&\quad + C_\epsilon\psi'F^{kk}|u_k|^2 + \epsilon C\psi'\mathcal{F}. \end{aligned}
\end{equation}

\subsection*{$\bullet$ The terms involving $|g_{1\bar 1k}|$ and $F^{kk}u_{k\bar k}$}
From \eqref{eq:b10} we know that
\[ g_{1\bar 1k} = \ti{\lambda}_{1,k} + O(\lambda_1), \]
and so using \eqref{eq:b142} we get
\begin{equation} \label{eq:e1} F^{kk}\lambda_1^{-1}|g_{1\bar 1k}| \leq \frac{1}{2K^{1/2}} \sum_pF^{kk}(|u_{\bar pk}| + |u_{pk}|) -
C_\epsilon \psi' F^{kk}|u_k|^2 - \epsilon \psi'\mathcal{F} +
C\mathcal{F}. \end{equation}
The terms involving $|u_{\bar pk}|, |u_{pk}|$ can be absorbed by the
squared terms $|u_{p\bar k}|^2, |u_{pk}|^2$ in (\ref{eq:star}), and so we obtain
\begin{equation}
\begin{aligned} \label{eq:d2} 0 &\geq \sum_p\frac{F^{kk}}{20K}(|u_{pk}|^2 + |u_{p\bar k}|^2) +
\frac{1}{2}\psi'' F^{kk}|u_k|^2 + \psi' F^{kk}u_{k\bar k} \\
&\quad - 2\psi'^2\delta^{-1}F^{11}K - C\mathcal{F}
+ C_\epsilon\psi'F^{kk}|u_k|^2 + \epsilon C\psi'\mathcal{F}. \end{aligned}\end{equation}
As for the term involving $u_{k\bar k}$, we have
  \[ \begin{aligned}\psi' F^{kk}u_{k\bar k} &= \psi' F^{kk}(g_{k\bar k} - \chi_{k\bar
    k} - W_{k\bar k}). \end{aligned}\]
  As in \eqref{eq:d1} we have
  \[ \sum_k F^{kk}W_{k\bar k} = \sum_i \widetilde{F}^{ii} Z_{i\bar
    i}. \]
Recall that $Z_{1\bar 1}$ does not contain $u_1$ or $u_{\bar
    1}$ and $\widetilde{F}^{11}$ is the only ``large'' coefficient, of order $F^{kk}$ for $k>1$. It
  follows that
  \[ |F^{kk} W_{k\bar k}| \leq CF^{kk} |u_k| \leq
  C_\epsilon F^{kk}|u_k|^2 + \epsilon
  \mathcal{F},\]
  and so
  \begin{equation}
  \begin{aligned}\label{eq:e2}
    \psi' F^{kk}u_{k\bar k} &\geq \psi' F^{kk}(g_{k\bar k} - \chi_{k\bar
    k}) + C_\epsilon\psi' F^{kk}|u_k|^2 + \epsilon\psi' \mathcal{F}. \end{aligned}\end{equation}
  From \eqref{eq:d2} we then finally obtain (if necessary replacing $C_\epsilon$ by another constant depending only on $\ve$ and the allowed data),
that
  \begin{equation}
   \begin{aligned}\label{eq:d3}
     0 &\geq  F^{11}\left(\frac{\lambda_1^2}{40K} -
    2\psi'^2\delta^{-1}K\right)  +
\left(\frac{1}{2}\psi'' + C_\epsilon\psi'\right) F^{kk}|u_k|^2  \\
&\quad - C_0 \mathcal{F}
+ \epsilon C_0 \psi'\mathcal{F} -\psi' F^{kk}(\chi_{k\bar k} - g_{k\bar k}),
 \end{aligned} \end{equation}
  for a uniform $C_0$.  We have used the fact that $|u_{1\ov{1}}|^2 \ge \frac{1}{2} \lambda_1^2 - CK$.

Under the assumption that the function $\underline{u}=0$ is a
$\mathcal{C}$-subsolution, and that $\lambda_1 \gg 1$, we may apply Proposition \ref{prop:subsol} and see that there is a
uniform positive number $\kappa > 0$ such that one of two possibilities occurs:
\begin{itemize}

\item[(a)] We have $F^{kk}(\chi_{k\bar k} - g_{k\bar k}) > \kappa \mathcal{F}$. In this
  case we have
  \[ \begin{aligned}0 &\geq F^{11}\left(\frac{\lambda_1^2}{40K} -
    2\psi'^2\delta^{-1}K\right) +
\left( \frac{1}{2}\psi'' + C_\epsilon\psi' \right) F^{kk}|u_k|^2 \\
&\quad - C_0\mathcal{F}
+ \epsilon C_0\psi'\mathcal{F} - \psi'\kappa\mathcal{F}. \end{aligned}\]
  We first choose $\epsilon > 0$ such that $\epsilon C_0 < \kappa / 2$. We then choose the parameter $D_2$ in the definition of $\psi(t)
  = D_1e^{-D_2t}$ to be large enough so that
  \[ \frac{1}{2}\psi'' > C_\epsilon|\psi'|. \]
  At this point we have
  \[ \begin{aligned}0 &\geq F^{11}\left(\frac{\lambda_1^2}{40K} -
    2\psi'^2\delta^{-1}K\right) - C_0\mathcal{F}
- \frac{1}{2}\psi'\kappa\mathcal{F}. \end{aligned}\]
  We now choose $D_1$ so large that $-\frac{1}{2} \psi'\kappa > C_0$, which
  implies
  \[ \frac{\lambda_1^2}{40K} \leq 2\psi'^2\delta^{-1}K. \]
  Note that $\delta$ is determined by the choices of $D_1,D_2$, according to \eqref{eq:definedelta}, so we
  obtain the required upper bound for $\lambda_1 / K$.

\medskip

  \item[(b)] We have $F^{11} > \kappa
  \mathcal{F}$. With the choices of constants made above,
  \eqref{eq:d3} implies that
  \[ \begin{aligned}0 &\geq \kappa\mathcal{F}\left(\frac{\lambda_1^2}{40K} -
    2\psi'^2\delta^{-1}K\right) - C_0\mathcal{F} \\
&\quad\quad +
  \epsilon C_0\psi'\mathcal{F} + C_1\psi'\mathcal{F} + \psi'F^{kk}g_{k\bar
    k}, \end{aligned}\]
    for another uniform constant $C_1$.
  Since $F^{kk}g_{k\bar k} \leq \mathcal{F}\lambda_1$, we can divide
  through by $\mathcal{F}K$ and obtain
  \[ 0\geq \frac{\kappa\lambda_1^2}{40K^2} - C_2(1 + K^{-1} +
  \lambda_1K^{-1}), \]
  for a uniform $C_2$.
  The required upper bound for $\lambda_1 / K$ follows from this.

\end{itemize}

\bigskip\noindent{\bf Case 2.} We now assume that $\delta \lambda_1 <
-\lambda_n$, with all the constants $D_1,D_2,\delta$ fixed as in the
previous case. We first use that  $F^{nn} \geq \frac{1}{n}\mathcal{F}$ as well as
$\lambda_n^2 > \delta^2\lambda_1^2$ to bound
\[\begin{aligned}
\sum_p\frac{F^{kk}}{6K}(|u_{pk}|^2
  + |u_{p\bar k}|^2)&\geq \frac{F^{nn}}{6K}|u_{n\ov{n}}|^2\geq \frac{\mathcal{F}}{6nK}|\lambda_n-\chi_{n\ov{n}}-W_{n\ov{n}}|^2\\
  &\geq \frac{\mathcal{F}}{10nK}|\lambda_n|^2-\frac{C\mathcal{F}}{K}(1+K)\\
  &\geq \frac{\delta^2}{10nK}\mathcal{F}\lambda_1^2-C\mathcal{F}.
  \end{aligned}
\]
In  \eqref{eq:b7}  we now discard the positive first
term and the term involving $\psi''$, and use this to obtain
\[\begin{aligned} 0 \geq {} & - \frac{F^{kk}|\ti{\lambda}_{1,k}|^2}{\lambda_1^2} +
\frac{\delta^2}{10nK}\mathcal{F}\lambda_1^2
 + \phi'' F^{kk}|V_k|^2 \\  & {} +
\psi'F^{kk}u_{k\bar k} - C\Big[F^{kk}\lambda_1^{-1}|g_{1\bar 1k}| +
\mathcal{F}\Big]. \end{aligned}\]

To deal with the terms involving $F^{kk}u_{k\bar k}$
 and $|g_{1\bar 1k}|$ we note that
 \[ F^{kk} |u_{k\bar k}| \le C \mathcal{F} \lambda_1 \]
 and, since $g_{1\ov{1}k} = \tilde{\lambda}_{1,k} + O(\lambda_1)$,
 \[ CF^{kk} \lambda_1^{-1} |g_{1\ov{1}k}| \leq CF^{kk} \lambda_1^{-1} |\tilde{\lambda}_{1,k}| + C \mathcal{F} \leq
   \frac{1}{2} \frac{F^{kk} |\tilde{\lambda}_{1,k}|^2}{\lambda_1^2} + C\mathcal{F}.\]
Then we obtain
\begin{equation}
 \begin{aligned}\label{eq:e4} 0 &\geq -\frac{3}{2}\frac{F^{kk}|\ti{\lambda}_{1,k}|^2}{\lambda_1^2} +
\frac{\delta^2}{10nK} \mathcal{F}\lambda_1^2 + F^{kk}\phi''|V_k|^2   - C \mathcal{F} \lambda_1 \end{aligned}\end{equation}
  Using $\ti{H}_k=0$ we have, since $\psi'$ is fixed now and bounded,
  \[\begin{aligned} \frac{3}{2}\frac{F^{kk}|\ti{\lambda}_{1,k}|^2}{\lambda_1^2} &= \frac{3}{2} F^{kk}| \phi'V_k +
  \psi'u_k|^2 \\
  &\leq 2F^{kk}\phi'^2|V_k|^2 + CF^{kk}\psi'^2|u_k|^2 \\
  &\leq F^{kk}\phi''|V_k|^2 + C\mathcal{F}K. \end{aligned}\]
Returning to  \eqref{eq:e4}, we obtain, since we may assume $\lambda_1\ge K$,
  \[ 0 \geq \frac{\delta^2\lambda_1^2}{10nK}\mathcal{F} -
  C\lambda_1 \mathcal{F}. \]
Dividing by $\lambda_1 \mathcal{F}$ gives the required bound for $\lambda_1 / K$.

Then we immediately deduce  the bound (\ref{mainbd}), namely
\begin{equation}
\label{bound}\sup_M |\ddbar u|_{\alpha}\leq C(\sup_M|\nabla u|^2_\alpha+1).\end{equation}
A blow-up argument as in \cite[Section 6]{Sz} combined with a Liouville theorem \cite[Section 5]{Sz} (see also \cite{DK, TW2,TW3}), shows that $\sup_M|\nabla u|^2_\alpha\leq C$
and so we get a uniform bound $|\Delta u|\leq C$. Here we remark that in the blow-up argument the only difference from the setup here (compared to \cite{Sz}) is the presence of the term $W_{i\ov{j}}$. However this term is linear in $\nabla u$ and so converges to zero uniformly on compact sets under the rescaling procedure of \cite{Sz} (compare \cite[Section 6]{TW3}).

We can then apply the Evans-Krylov-type result in \cite[Theorem 1.1]{TWWY} and deduce a uniform bound
$$\|u\|_{C^{2,\beta}(M,\alpha)}\leq C,$$ for a uniform $0<\beta<1$. Differentiating the equation and applying a standard bootstrapping argument we finally obtain uniform higher-order estimates.

\section{Proof of Theorem \ref{main2}}  \label{sec:last}

In this section, we explain how Theorem \ref{main2} follows from  Theorem \ref{main}.

Write $*$ for the Hodge star operator with respect to $\alpha$.   This acts on real $(n-1, n-1)$ forms as follows.  Consider a real $(n-1, n-1)$ form $\Theta$ given by
\[\begin{aligned}
\Theta = {} & (\sqrt{-1})^{n-1} \sum_{i,j} (\textrm{sgn}(i,j)) \Theta_{i\ov{j}} dz^1 \wedge d\ov{z}^1 \wedge \cdots \wedge \widehat{dz^i} \wedge d\ov{z}^i \wedge \cdots \\
& {} \wedge dz^j \wedge \widehat{d\ov{z}^j} \wedge \cdots \wedge dz^n \wedge d\ov{z}^n,
\end{aligned}
\]
with $\textrm{sgn}(i,j)=1$ for $i\le j$ and $\textrm{sgn}(i,j)=-1$ if $i > j$.  If we are computing at a point in coordinates so that $\alpha_{i\ov{j}}=\delta_{ij}$, then
$$*\Theta = \sqrt{-1} \sum_{i,j} \Theta_{i\ov{j}} dz^j \wedge d\ov{z}^i.$$
A basic property is that for any Hermitian metric $\omega$ we have (see \cite[Section 2]{TW2}, for example)
$$\left( \frac{\omega^n}{\alpha^n} \right)^{n-1} = \frac{ (* (\omega^{n-1}))^n}{(*(\alpha^{n-1}))^n}= \frac{ (* (\omega^{n-1}))^n}{((n-1)!\alpha)^n}.$$
Then taking $\omega$ as in Theorem \ref{main2},  we see that equation (\ref{cma}) is equivalent to
\begin{equation}
 \label{eq:eq}
\log \frac{ (*(\omega^{n-1}))^n}{((n-1)!\alpha)^n} = h,
\end{equation}
with  $h= (n-1)(F+b)$, a smooth function.  Recall that
$$\omega^{n-1}=\alpha_0^{n-1}+\ddbar u \wedge\alpha^{n-2}+\mathrm{Re}\left(\mn \de u \wedge\db(\alpha^{n-2})\right).$$
As in \cite{TW2}, we have
$$\frac{1}{(n-1)!}*(\ddbar u \wedge \alpha^{n-2}) = \frac{1}{n-1} ((\Delta u)\alpha - \ddbar u).$$
Define
\begin{equation}
 \label{eq:Zdefn}
Z_{i\ov{j}} = \left(\frac{1}{(n-1)!}*\mathrm{Re}\left(\mn \de u \wedge\db(\alpha^{n-2})\right) \right)_{i\ov{j}}.
\end{equation}
A straightforward but long calculation gives
\begin{equation}\label{eq:Zdefn2}
\begin{split}
Z_{i\ov{j}} =  {}   \frac{1}{2(n-1)} \bigg\{ &\alpha^{p\ov{q}} \alpha^{k\ov{\ell}} u_p \ov{T_{q \ell \ov{k}}} \alpha_{i\ov{j}} - \alpha^{k\ov{\ell}} u_{i} \ov{T_{j\ell \ov{k}}} - \alpha^{k\ov{\ell}} u_{k} \ov{T_{\ell j \ov{i}}}  \\
&  + \alpha^{p\ov{q}} \alpha^{k\ov{\ell}} u_{\ov{q}} T_{pk\ov{\ell}} \alpha_{i\ov{j}} - \alpha^{k\ov{\ell}} u_{\ov{j}} T_{ik\ov{\ell}} -\alpha^{k\ov{\ell}} u_{\ov{\ell}} T_{k i\ov{j}} \bigg\},
\end{split}
\end{equation}
where we are writing $T_{ij}^k$ for the torsion of $\alpha$ and $T_{ij\ov{\ell}}=T_{ij}^k \alpha_{k\ell}$.
An important point to note is that, since the torsion is skew-symmetric $T_{ij\ov{\ell}}=-T_{ji\ov{\ell}}$, in orthonormal coordinates for $\alpha$ we see that $Z_{i\ov{j}}$ is independent of $u_{\ov{i}}$ and $u_j$, and that $\nabla_i Z_{i\ov{i}}$ is independent of $u_{\ov{i}}$. Indeed, in local orthonormal coordinates for $\alpha$ we have
$$Z_{i\ov{i}}=\frac{1}{2(n-1)} \left(\sum_{p\neq i}\sum_{k\neq i}u_p \ov{T_{pk\ov{k}}}+\sum_{p\neq i}\sum_{k\neq i}u_{\ov{p}} T_{pk\ov{k}}\right),$$
and for $i\neq j$
$$Z_{i\ov{j}}=\frac{1}{2(n-1)} \left(-\sum_{k\neq j}(u_i \ov{T_{jk\ov{k}}}+u_k \ov{T_{kj\ov{i}}})-\sum_{k\neq i}(u_{\ov{j}} T_{ik\ov{k}}+u_{\ov{k}} T_{ki\ov{j}})\right),$$
using the skew-symmetry of the torsion. Also,
\[\begin{split}
\nabla_i Z_{i\ov{i}}=\frac{1}{2(n-1)} \bigg(
&\sum_{p\neq i}\sum_{k\neq i}( u_p \nabla_i\ov{T_{pk\ov{k}}}+\nabla_i u_p \ov{T_{pk\ov{k}}})\\
&+\sum_{p\neq i}\sum_{k\neq i}(u_{\ov{p}} \nabla_iT_{pk\ov{k}} +\nabla_i u_{\ov{p}} T_{pk\ov{k}} )\bigg),
\end{split}\]
and the statement follows.
 We also define
$$\tilde{\chi}_{i\ov{j}} = \left(\frac{1}{(n-1)!} *(\alpha_0^{n-1})\right)_{i\ov{j}}.$$
Given this, we see that \eqref{eq:eq} is equivalent to
$$\log ( \mu_1 \cdots \mu_n) = h,$$
where $\mu_i$ are the eigenvalues of $\alpha^{i\ov{p}} \tilde{g}_{j\ov{p}}$, for $\tilde{g}$ given by
$$\tilde{g}_{i\ov{j}} = \tilde{\chi}_{i\ov{j}} + \frac{1}{n-1} ((\Delta u)\alpha_{i\ov{j}} - u_{i\ov{j}} ) +Z_{i\ov{j}}.$$
Since $\ti{\chi}_{i\ov{j}}$ is positive definite, we have that $0$ is a $\mathcal{C}$-subsolution.
From the discussion in section \ref{sec:bg}, it is now immediate to see that this equation falls into the setup of Theorem \ref{main}, and so we obtain  the uniform a priori estimate \eqref{bound}. Therefore Theorem \ref{main2} follows from \cite[Theorem 1.7]{TW3}.

\end{document}